\documentclass[12pt,a4paper]{article}

\usepackage{amsfonts,amssymb,amsthm,amsmath,latexsym}
\usepackage{subeqnarray,eqsection,indent,url,cite,bm}
\usepackage{multirow}  
\usepackage{tensor}  
\usepackage{float}  
\usepackage{graphicx} 

\usepackage{calc}
\usepackage{pgf,tikz}

\date{January 6, 2020 \\[1.5mm] revised March 10, 2020}   

\begin{document}

\title{\vspace*{-2cm}Phylogenetic trees, augmented perfect matchings, \\
       and a Thron-type continued fraction (T-fraction) \\
       for the Ward polynomials}

\author{
     \\
     {\small Andrew Elvey Price}             \\[-2mm]
     {\small\it LaBRI}                       \\[-2mm]
     {\small\it Universit\'e de Bordeaux}     \\[-2mm]
     {\small\it F-33405 Talence Cedex}        \\[-2mm]
     {\small\it FRANCE}                   \\[-2mm]
     {\small\tt andrewelveyprice@gmail.com}  \\[-2mm]
     {\protect\makebox[5in]{\quad}}  
     \\[2mm]
     {\small Alan D.~Sokal}                  \\[-2mm]
     {\small\it Department of Mathematics}   \\[-2mm]
     {\small\it University College London}   \\[-2mm]
     {\small\it Gower Street}                \\[-2mm]
     {\small\it London WC1E 6BT}             \\[-2mm]
     {\small\it UNITED KINGDOM}              \\[-2mm]
     {\small\tt sokal@math.ucl.ac.uk}        \\[-2mm]
     \\[-2mm]
     {\small\it Department of Physics}       \\[-2mm]
     {\small\it New York University}         \\[-2mm]
     {\small\it 726 Broadway}          \\[-2mm]
     {\small\it New York, NY 10003}      \\[-2mm]
     {\small\it USA}      \\[-2mm]
     {\small\tt sokal@nyu.edu}            \\[3mm]
}

\maketitle
\thispagestyle{empty}   

\begin{abstract}
We find a Thron-type continued fraction (T-fraction)
for the ordinary generating function of the Ward polynomials,
as well as for some generalizations employing a large (indeed infinite)
family of independent indeterminates.
Our proof is based on a bijection between
super-augmented perfect matchings and labeled Schr\"oder paths,
which generalizes Flajolet's bijection between
perfect matchings and labeled Dyck paths.
\end{abstract}


\bigskip
\noindent
{\bf Mathematics Subject Classification (MSC 2010) codes:}
05A19 (Primary);
05A10, 05A15, 05A18, 30B70, 92B10 (Secondary).

\clearpage

\newtheorem{theorem}{Theorem}[section]
\newtheorem{proposition}[theorem]{Proposition}
\newtheorem{lemma}[theorem]{Lemma}
\newtheorem{corollary}[theorem]{Corollary}
\newtheorem{definition}[theorem]{Definition}
\newtheorem{conjecture}[theorem]{Conjecture}
\newtheorem{question}[theorem]{Question}
\newtheorem{problem}[theorem]{Problem}
\newtheorem{example}[theorem]{Example}

\renewcommand{\theenumi}{\alph{enumi}}
\renewcommand{\labelenumi}{(\theenumi)}
\def\eop{\hbox{\kern1pt\vrule height6pt width4pt
depth1pt\kern1pt}\medskip}
\def\prf{\par\noindent{\bf Proof.\enspace}\rm}
\def\rmk{\par\medskip\noindent{\bf Remark\enspace}\rm}

\newcommand{\bigdash}{%
\smallskip\begin{center} \rule{5cm}{0.1mm} \end{center}\smallskip}

\newcommand{\be}{\begin{equation}}
\newcommand{\ee}{\end{equation}}
\newcommand{\<}{\langle}
\renewcommand{\>}{\rangle}
\newcommand{\widebar}{\overline}
\def\reff#1{(\protect\ref{#1})}
\def\spose#1{\hbox to 0pt{#1\hss}}
\def\ltapprox{\mathrel{\spose{\lower 3pt\hbox{$\mathchar"218$}}
    \raise 2.0pt\hbox{$\mathchar"13C$}}}
\def\gtapprox{\mathrel{\spose{\lower 3pt\hbox{$\mathchar"218$}}
    \raise 2.0pt\hbox{$\mathchar"13E$}}}
\def\textprime{${}^\prime$}
\def\proof{\par\medskip\noindent{\sc Proof.\ }}
\def\firstproof{\par\medskip\noindent{\sc First Proof.\ }}
\def\secondproof{\par\medskip\noindent{\sc Second Proof.\ }}
\def\qed{ $\square$ \bigskip}
\newcommand{\myendremark}{ $\blacksquare$ \bigskip}
\def\proofof#1{\bigskip\noindent{\sc Proof of #1.\ }}
\def\firstproofof#1{\bigskip\noindent{\sc First Proof of #1.\ }}
\def\secondproofof#1{\bigskip\noindent{\sc Second Proof of #1.\ }}
\def\thirdproofof#1{\bigskip\noindent{\sc Third Proof of #1.\ }}
\def\half{ {1 \over 2} }
\def\third{ {1 \over 3} }
\def\twothird{ {2 \over 3} }
\def\smfrac#1#2{{\textstyle{#1\over #2}}}
\def\smhalf{ {\smfrac{1}{2}} }
\newcommand{\real}{\mathop{\rm Re}\nolimits}
\renewcommand{\Re}{\mathop{\rm Re}\nolimits}
\newcommand{\imag}{\mathop{\rm Im}\nolimits}
\renewcommand{\Im}{\mathop{\rm Im}\nolimits}
\newcommand{\sgn}{\mathop{\rm sgn}\nolimits}
\newcommand{\tr}{\mathop{\rm tr}\nolimits}
\newcommand{\supp}{\mathop{\rm supp}\nolimits}
\newcommand{\disc}{\mathop{\rm disc}\nolimits}
\newcommand{\diag}{\mathop{\rm diag}\nolimits}
\newcommand{\tridiag}{\mathop{\rm tridiag}\nolimits}
\newcommand{\AZ}{\mathop{\rm AZ}\nolimits}
\newcommand{\perm}{\mathop{\rm perm}\nolimits}
\def\hboxscript#1{ {\hbox{\scriptsize\em #1}} }
\renewcommand{\emptyset}{\varnothing}
\newcommand{\eqdef}{\stackrel{\rm def}{=}}

\newcommand{\restrict}{\!\upharpoonright\!}

\newcommand{\compinv}{{\langle -1 \rangle}}   

\newcommand{\scra}{{\mathcal{A}}}
\newcommand{\scrb}{{\mathcal{B}}}
\newcommand{\scrc}{{\mathcal{C}}}
\newcommand{\scrd}{{\mathcal{D}}}
\newcommand{\scre}{{\mathcal{E}}}
\newcommand{\scrf}{{\mathcal{F}}}
\newcommand{\scrg}{{\mathcal{G}}}
\newcommand{\scrh}{{\mathcal{H}}}
\newcommand{\scri}{{\mathcal{I}}}
\newcommand{\scrk}{{\mathcal{K}}}
\newcommand{\scrl}{{\mathcal{L}}}
\newcommand{\scrm}{{\mathcal{M}}}
\newcommand{\scrn}{{\mathcal{N}}}
\newcommand{\scro}{{\mathcal{O}}}
\newcommand{\scrp}{{\mathcal{P}}}
\newcommand{\scrq}{{\mathcal{Q}}}
\newcommand{\scrr}{{\mathcal{R}}}
\newcommand{\scrs}{{\mathcal{S}}}
\newcommand{\scrt}{{\mathcal{T}}}
\newcommand{\scrv}{{\mathcal{V}}}
\newcommand{\scrw}{{\mathcal{W}}}
\newcommand{\scrz}{{\mathcal{Z}}}

\newcommand{\ahat}{{\widehat{a}}}
\newcommand{\Zhat}{{\widehat{Z}}}
\renewcommand{\k}{{\mathbf{k}}}
\newcommand{\n}{{\mathbf{n}}}
\newcommand{\vv}{{\mathbf{v}}}
\newcommand{\bv}{{\mathbf{v}}}
\newcommand{\w}{{\mathbf{w}}}
\newcommand{\x}{{\mathbf{x}}}
\newcommand{\bz}{{\mathbf{z}}}
\newcommand{\bw}{{\mathbf{w}}}
\newcommand{\cc}{{\mathbf{c}}}
\newcommand{\zero}{{\mathbf{0}}}
\newcommand{\one}{{\mathbf{1}}}
\newcommand{\bmm}{ {\bf m} }

\newcommand{\C}{{\mathbb C}}
\newcommand{\D}{{\mathbb D}}
\newcommand{\Z}{{\mathbb Z}}
\newcommand{\N}{{\mathbb N}}
\newcommand{\Q}{{\mathbb Q}}
\newcommand{\PP}{{\mathbb P}}
\newcommand{\R}{{\mathbb R}}
\newcommand{\RR}{{\mathbb R}}
\newcommand{\E}{{\mathbb E}}

\newcommand{\Sym}{{\mathfrak{S}}}
\newcommand{\SymB}{{\mathfrak{B}}}

\newcommand{\myle}{\preceq}
\newcommand{\myge}{\succeq}
\newcommand{\mygt}{\succ}

\newcommand{\B}{{\sf B}}
\newcommand{\OB}{{\sf OB}}
\newcommand{\OS}{{\sf OS}}
\newcommand{\OO}{{\sf O}}
\newcommand{\SP}{{\sf SP}}
\newcommand{\OSP}{{\sf OSP}}
\newcommand{\Eu}{{\sf Eu}}
\newcommand{\ERR}{{\sf ERR}}
\newcommand{\sfB}{{\sf B}}
\newcommand{\sfD}{{\sf D}}
\newcommand{\sfE}{{\sf E}}
\newcommand{\sfG}{{\sf G}}
\newcommand{\sfJ}{{\sf J}}
\newcommand{\sfP}{{\sf P}}
\newcommand{\sfQ}{{\sf Q}}
\newcommand{\sfS}{{\sf S}}
\newcommand{\sfT}{{\sf T}}
\newcommand{\sfW}{{\sf W}}
\newcommand{\sfMV}{{\sf MV}}
\newcommand{\AMV}{{\sf AMV}}
\newcommand{\BM}{{\sf BM}}

\newcommand{\emIB}{{\hbox{\em IB}}}
\newcommand{\emIP}{{\hbox{\em IP}}}
\newcommand{\emOB}{{\hbox{\em OB}}}
\newcommand{\emSC}{{\hbox{\em SC}}}

\newcommand{\clop}{{\rm clop}}
\newcommand{\stat}{{\rm stat}}
\newcommand{\cyc}{{\rm cyc}}
\newcommand{\Asc}{{\rm Asc}}
\newcommand{\asc}{{\rm asc}}
\newcommand{\Des}{{\rm Des}}
\newcommand{\des}{{\rm des}}
\newcommand{\Exc}{{\rm Exc}}
\newcommand{\exc}{{\rm exc}}
\newcommand{\aexc}{{\rm aexc}}
\newcommand{\Wex}{{\rm Wex}}
\newcommand{\wex}{{\rm wex}}
\newcommand{\Fix}{{\rm Fix}}
\newcommand{\fix}{{\rm fix}}
\newcommand{\bfix}{{\mathbf{fix}}}
\newcommand{\lev}{{\rm lev}}
\newcommand{\lrmax}{{\rm lrmax}}
\newcommand{\rlmax}{{\rm rlmax}}
\newcommand{\Rec}{{\rm Rec}}
\newcommand{\rec}{{\rm rec}}
\newcommand{\Arec}{{\rm Arec}}
\newcommand{\arec}{{\rm arec}}
\newcommand{\ERec}{{\rm ERec}}
\newcommand{\erec}{{\rm erec}}
\newcommand{\EArec}{{\rm EArec}}
\newcommand{\earec}{{\rm earec}}
\newcommand{\recarec}{{\rm recarec}}
\newcommand{\nonrec}{{\rm nonrec}}
\newcommand{\nrar}{{\rm nrar}}
\newcommand{\ereccval}{{\rm ereccval}}
\newcommand{\ereccdrise}{{\rm ereccdrise}}
\newcommand{\eareccpeak}{{\rm eareccpeak}}
\newcommand{\eareccdfall}{{\rm eareccdfall}}
\newcommand{\rar}{{\rm rar}}
\newcommand{\nrcpeak}{{\rm nrcpeak}}
\newcommand{\nrcval}{{\rm nrcval}}
\newcommand{\nrcdrise}{{\rm nrcdrise}}
\newcommand{\nrcdfall}{{\rm nrcdfall}}
\newcommand{\nrfix}{{\rm nrfix}}
\newcommand{\Cpeak}{{\rm Cpeak}}
\newcommand{\cpeak}{{\rm cpeak}}
\newcommand{\Cval}{{\rm Cval}}
\newcommand{\cval}{{\rm cval}}
\newcommand{\Cdasc}{{\rm Cdasc}}
\newcommand{\cdasc}{{\rm cdasc}}
\newcommand{\Cddes}{{\rm Cddes}}
\newcommand{\cddes}{{\rm cddes}}
\newcommand{\cdrise}{{\rm cdrise}}
\newcommand{\cdfall}{{\rm cdfall}}
\newcommand{\cross}{{\rm cross}}
\newcommand{\nest}{{\rm nest}}
\newcommand{\ucross}{{\rm ucross}}
\newcommand{\ucrosscval}{{\rm ucrosscval}}
\newcommand{\ucrosscdrise}{{\rm ucrosscdrise}}
\newcommand{\lcross}{{\rm lcross}}
\newcommand{\lcrosscpeak}{{\rm lcrosscpeak}}
\newcommand{\lcrosscdfall}{{\rm lcrosscdfall}}
\newcommand{\unest}{{\rm unest}}
\newcommand{\unestcval}{{\rm unestcval}}
\newcommand{\unestcdrise}{{\rm unestcdrise}}
\newcommand{\lnest}{{\rm lnest}}
\newcommand{\lnestcpeak}{{\rm lnestcpeak}}
\newcommand{\lnestcdfall}{{\rm lnestcdfall}}
\newcommand{\ujoin}{{\rm ujoin}}
\newcommand{\ljoin}{{\rm ljoin}}
\newcommand{\upsnest}{{\rm upsnest}}
\newcommand{\lpsnest}{{\rm lpsnest}}
\newcommand{\psnest}{{\rm psnest}}
\newcommand{\ecp}{{\rm ecp}}
\newcommand{\ecar}{{\rm ecar}}
\newcommand{\car}{{\rm car}}
\newcommand{\ecnar}{{\rm ecnar}}
\newcommand{\cnar}{{\rm cnar}}
\newcommand{\ocp}{{\rm ocp}}
\newcommand{\ocar}{{\rm ocar}}
\newcommand{\ocnar}{{\rm ocnar}}
\newcommand{\eor}{{\rm eor}}
\newcommand{\eonr}{{\rm eonr}}
\newcommand{\oor}{{\rm oor}}
\newcommand{\oonr}{{\rm oonr}}
\newcommand{\wig}{{\rm wig}}
\newcommand{\dash}{{\rm dash}}
\newcommand{\Peak}{{\rm Peak}}
\newcommand{\peak}{{\rm peak}}
\newcommand{\Val}{{\rm Val}}
\newcommand{\val}{{\rm val}}
\newcommand{\Dasc}{{\rm Dasc}}
\newcommand{\dasc}{{\rm dasc}}
\newcommand{\Ddes}{{\rm Ddes}}
\newcommand{\ddes}{{\rm ddes}}
\newcommand{\inv}{{\rm inv}}
\newcommand{\maj}{{\rm maj}}
\newcommand{\rs}{{\rm rs}}
\newcommand{\crr}{{\rm cr}}
\newcommand{\crosshat}{{\widehat{\rm cr}}}
\newcommand{\nee}{{\rm ne}}
\newcommand{\qne}{{\rm qne}}
\newcommand{\psne}{{\rm psne}}
\newcommand{\crne}{{\rm crne}}
\newcommand{\ov}{{\rm ov}}
\newcommand{\cov}{{\rm cov}}
\newcommand{\pscov}{{\rm pscov}}
\newcommand{\rodd}{{\rm rodd}}
\newcommand{\reven}{{\rm reven}}
\newcommand{\lodd}{{\rm lodd}}
\newcommand{\leven}{{\rm leven}}
\newcommand{\sg}{{\rm sg}}
\newcommand{\bl}{{\rm bl}}
\newcommand{\tran}{{\rm tr}}
\newcommand{\area}{{\rm area}}
\newcommand{\ret}{{\rm ret}}
\newcommand{\peaks}{{\rm peaks}}
\newcommand{\hl}{{\rm hl}}
\newcommand{\sll}{{\rm sll}}
\newcommand{\negg}{{\rm neg}}
\newcommand{\imp}{{\rm imp}}

\newcommand{\sfa}{{{\sf a}}}
\newcommand{\sfb}{{{\sf b}}}
\newcommand{\sfc}{{{\sf c}}}
\newcommand{\sfd}{{{\sf d}}}
\newcommand{\sfe}{{{\sf e}}}
\newcommand{\sff}{{{\sf f}}}
\newcommand{\sfg}{{{\sf g}}}
\newcommand{\sfh}{{{\sf h}}}
\newcommand{\sfi}{{{\sf i}}}
\newcommand{\bsfa}{{\mbox{\textsf{\textbf{a}}}}}
\newcommand{\bsfb}{{\mbox{\textsf{\textbf{b}}}}}
\newcommand{\bsfc}{{\mbox{\textsf{\textbf{c}}}}}
\newcommand{\bsfd}{{\mbox{\textsf{\textbf{d}}}}}
\newcommand{\bsfe}{{\mbox{\textsf{\textbf{e}}}}}
\newcommand{\bsff}{{\mbox{\textsf{\textbf{f}}}}}
\newcommand{\bsfg}{{\mbox{\textsf{\textbf{g}}}}}
\newcommand{\bsfh}{{\mbox{\textsf{\textbf{h}}}}}
\newcommand{\bsfi}{{\mbox{\textsf{\textbf{i}}}}}

\newcommand{\bfx}{{\bf x}}

\newcommand{\ba}{{\bm{a}}}
\newcommand{\bahat}{{\widehat{\bm{a}}}}
\newcommand{\bb}{{\bm{b}}}
\newcommand{\bc}{{\bm{c}}}
\newcommand{\bff}{{\bm{f}}}
\newcommand{\bg}{{\bm{g}}}
\newcommand{\br}{{\bm{r}}}
\newcommand{\bu}{{\bm{u}}}
\newcommand{\bA}{{\bm{A}}}
\newcommand{\bB}{{\bm{B}}}
\newcommand{\bC}{{\bm{C}}}
\newcommand{\bE}{{\bm{E}}}
\newcommand{\bF}{{\bm{F}}}
\newcommand{\bI}{{\bm{I}}}
\newcommand{\bJ}{{\bm{J}}}
\newcommand{\bM}{{\bm{M}}}
\newcommand{\bN}{{\bm{N}}}
\newcommand{\bP}{{\bm{P}}}
\newcommand{\bQ}{{\bm{Q}}}
\newcommand{\bS}{{\bm{S}}}
\newcommand{\bT}{{\bm{T}}}
\newcommand{\bW}{{\bm{W}}}
\newcommand{\bIB}{{\bm{IB}}}
\newcommand{\bOB}{{\bm{OB}}}
\newcommand{\bOS}{{\bm{OS}}}
\newcommand{\bERR}{{\bm{ERR}}}
\newcommand{\bSP}{{\bm{SP}}}
\newcommand{\bMV}{{\bm{MV}}}
\newcommand{\bBM}{{\bm{BM}}}
\newcommand{\balpha}{{\bm{\alpha}}}
\newcommand{\bbeta}{{\bm{\beta}}}
\newcommand{\bgamma}{{\bm{\gamma}}}
\newcommand{\bdelta}{{\bm{\delta}}}
\newcommand{\bomega}{{\bm{\omega}}}
\newcommand{\bone}{{\bm{1}}}
\newcommand{\bzero}{{\bm{0}}}

\newcommand{\Cbar}{{\overline{C}}}
\newcommand{\Dbar}{{\overline{D}}}
\newcommand{\dbar}{{\overline{d}}}
\def\Ctilde{{\widetilde{C}}}
\def\Ftilde{{\widetilde{F}}}
\def\Gtilde{{\widetilde{G}}}
\def\Htilde{{\widetilde{H}}}
\def\Ptilde{{\widetilde{P}}}
\def\Chat{{\widehat{C}}}
\def\ctilde{{\widetilde{c}}}
\def\zbar{{\overline{Z}}}
\def\pitilde{{\widetilde{\pi}}}

\newcommand{\sech}{{\rm sech}}

%
%
\newcommand{\sn}{{\rm sn}}
\newcommand{\cn}{{\rm cn}}
\newcommand{\dn}{{\rm dn}}
\newcommand{\sm}{{\rm sm}}
\newcommand{\cm}{{\rm cm}}

%
%
\newcommand{\zfz}{ {{}_0 \! F_0} }
\newcommand{\zfo}{ {{}_0 \! F_1} }
\newcommand{\ofz}{ {{}_1 \! F_0} }
\newcommand{\oft}{ {{}_1 \! F_2} }

%
%
\newcommand{\FHyper}[2]{ {\tensor[_{#1 \!}]{F}{_{#2}}\!} }
\newcommand{\phiHyper}[2]{ {\tensor[_{#1}]{\phi}{_{#2}}} }
\newcommand{\FHYPER}[5]{ {\FHyper{#1}{#2} \!\biggl(
   \!\!\begin{array}{c} #3 \\[1mm] #4 \end{array}\! \bigg|\, #5 \! \biggr)} }
\newcommand{\ofo}{ {\FHyper{1}{1}} }
\newcommand{\tfo}{ {\FHyper{2}{1}} }
\newcommand{\FHYPERbottomzero}[3]{ {\FHyper{#1}{0} \!\biggl(
   \!\!\begin{array}{c} #2 \\[1mm] \hbox{---} \end{array}\! \bigg|\, #3 \! \biggr)} }
\newcommand{\FHYPERtopzero}[3]{ {\FHyper{0}{#1} \!\biggl(
   \!\!\begin{array}{c} \hbox{---} \\[1mm] #2 \end{array}\! \bigg|\, #3 \! \biggr)} }

%
%
\newcommand{\stirlingsubset}[2]{\genfrac{\{}{\}}{0pt}{}{#1}{#2}}
\newcommand{\stirlingcycleold}[2]{\genfrac{[}{]}{0pt}{}{#1}{#2}}
\newcommand{\stirlingcycle}[2]{\genfrac{[}{]}{0pt}{}{#1}{#2}}
\newcommand{\assocstirlingsubset}[3]{{\genfrac{\{}{\}}{0pt}{}{#1}{#2}}_{\! \ge #3}}
\newcommand{\assocstirlingsubsetin}[3]{{\genfrac{\{}{\}}{0pt}{}{#1}{#2}}_{\! \in #3}}
\newcommand{\genstirlingsubset}[4]{{\genfrac{\{}{\}}{0pt}{}{#1}{#2}}_{\! #3,#4}}
\newcommand{\euler}[2]{\genfrac{\langle}{\rangle}{0pt}{}{#1}{#2}}
\newcommand{\eulergen}[3]{{\genfrac{\langle}{\rangle}{0pt}{}{#1}{#2}}_{\! #3}}
\newcommand{\eulersecond}[2]{\left\langle\!\!\mathchoice{\!}{}{}{} \euler{#1}{#2} \mathchoice{\!}{}{}{}\!\!\right\rangle}
\newcommand{\eulersecondgen}[3]{{\left\langle\!\! \euler{#1}{#2} \!\!\right\rangle}_{\! #3}}
\newcommand{\binomvert}[2]{\genfrac{\vert}{\vert}{0pt}{}{#1}{#2}}
\newcommand{\binomsquare}[2]{\genfrac{[}{]}{0pt}{}{#1}{#2}}


\newenvironment{sarray}{
             \textfont0=\scriptfont0
             \scriptfont0=\scriptscriptfont0
             \textfont1=\scriptfont1
             \scriptfont1=\scriptscriptfont1
             \textfont2=\scriptfont2
             \scriptfont2=\scriptscriptfont2
             \textfont3=\scriptfont3
             \scriptfont3=\scriptscriptfont3
           \renewcommand{\arraystretch}{0.7}
           \begin{array}{l}}{\end{array}}

\newenvironment{scarray}{
             \textfont0=\scriptfont0
             \scriptfont0=\scriptscriptfont0
             \textfont1=\scriptfont1
             \scriptfont1=\scriptscriptfont1
             \textfont2=\scriptfont2
             \scriptfont2=\scriptscriptfont2
             \textfont3=\scriptfont3
             \scriptfont3=\scriptscriptfont3
           \renewcommand{\arraystretch}{0.7}
           \begin{array}{c}}{\end{array}}

\newcommand*{\Scale}[2][4]{\scalebox{#1}{$#2$}}

\usetikzlibrary{decorations.pathmorphing}
\tikzset{snake it/.style={decorate, decoration={snake,segment length=5,amplitude=1}}}

\newcommand{\plotpt}[3][] 
{ \fill[#1,radius=0.175] (#2,#3) circle; }

\newcommand{\plotpermnobox}[3][]  
{
  \foreach \y [count=\x] in {#3}
  {
    \ifnum0=\y {} \else {
      \plotpt[#1]{\x}{\y}
    } \fi
  }
}

\newcommand{\PartArc}[4] {\draw [black,very thick] (#2,1) arc ((360*(#3)):(360*(#4)):{0.5*((#2)-(#1))} and {0.33*((#2)-(#1))+0.2});}

\newcommand{\Arc}[2] {\draw [black,very thick] (#2,1) arc (0:180:{0.5*((#2)-(#1))} and {0.33*((#2)-(#1))+0.2});} 

\newcommand{\Dashed}[2] {\draw [black,thick,dashed] (#1,1)--(#2,1);}

\newcommand{\Wiggly}[2] {\draw [black,thick,snake it] (#1,1)--(#2,1);}

\clearpage

\section{Introduction}

If $(a_n)_{n \ge 0}$ is a sequence of combinatorial numbers or polynomials
with $a_0 = 1$, it is often fruitful to seek to express its
ordinary generating function as a continued fraction.
The most commonly studied types of continued fractions
are Stieltjes-type (S-fractions),
\be
   \sum_{n=0}^\infty a_n t^n
   \;=\;
   \cfrac{1}{1 - \cfrac{\alpha_1 t}{1 - \cfrac{\alpha_2 t}{1 - \cdots}}}
   \label{def.Stype}
   \;\;,
\ee
and Jacobi-type (J-fractions),
\be
   \sum_{n=0}^\infty a_n t^n
   \;=\;
   \cfrac{1}{1 - \gamma_0 t - \cfrac{\beta_1 t^2}{1 - \gamma_1 t - \cfrac{\beta_2 t^2}{1 - \cdots}}}
   \label{def.Jtype}
   \;\;.
\ee
(Both sides of these expressions are to be interpreted as
formal power series in the indeterminate $t$.)
This line of investigation goes back at least to
Euler \cite{Euler_1760,Euler_1788},
but it gained impetus following Flajolet's \cite{Flajolet_80}
seminal discovery that any S-fraction (resp.\ J-fraction)
can be interpreted combinatorially as a generating function
for Dyck (resp.\ Motzkin) paths with suitable weights for each rise and fall
(resp.\ each rise, fall and level step).
There are now literally dozens of sequences $(a_n)_{n \ge 0}$
of combinatorial numbers or polynomials for which
a continued-fraction expansion of the type \reff{def.Stype} or \reff{def.Jtype}
is explicitly known.

A less commonly studied type of continued fraction is 
the Thron-type (T-fraction):
\be
   \sum_{n=0}^\infty a_n t^n
   \;=\;
   \cfrac{1}{1 - \delta_1 t - \cfrac{\alpha_1 t}{1 - \delta_2 t - \cfrac{\alpha_2 t}{1 - \delta_3 t -  \cfrac{\alpha_3 t}{1- \cdots}}}}
   \label{def.Ttype}
   \;\;.
\ee
Clearly the T-fractions are a generalization of the S-fractions,
and reduce to them when $\delta_i = 0$ for all $i$.
Besides the S-fractions,
there are four special types of T-fractions that are comparatively simple:
\begin{itemize}
   \item[1)]
Euler \cite[Chapter~18]{Euler_1748}
(see also \cite[pp.~17--18]{Wall_48})
showed that an arbitrary sequence of nonzero real numbers
(or finite sequence of nonzero real numbers followed by zeroes)
can be written as a special type of T-fraction,
namely, one in which
$\delta_1 = 0$ and $\delta_i = -\alpha_{i-1}$ for $i \ge 2$:
\be
   \cfrac{1}{1 - \cfrac{\alpha_1 t}{1 + \alpha_1 t - \cfrac{\alpha_2 t}{1 + \alpha_2 t -  \cfrac{\alpha_3 t}{1- \cdots}}}}
   \;\;=\;\;
   \sum_{n=0}^\infty \alpha_1 \alpha_2 \cdots \alpha_n \, t^n
   \;.
 \label{eq.euler_continued_fraction.0}
\ee
   \item[2)]
If $\delta_n = 0$ for all {\em even}\/ $n$,
then the T-fraction reduces to a J-fraction via the contraction formula
\cite{Sokal_totalpos}
\begin{subeqnarray}
   \gamma_0  & = &  \alpha_1 + \delta_1   \\
   \gamma_n  & = &  \alpha_{2n} + \alpha_{2n+1} + \delta_{2n+1}
                            \qquad\hbox{for $n \ge 1$} \\
   \beta_n  & = &  \alpha_{2n-1} \alpha_{2n}
 \label{eq.contraction_even.Ttype.coeffs}
\end{subeqnarray}
(which generalizes a well-known contraction formula
 \cite[p.~21]{Wall_48} \cite[p.~V-31]{Viennot_83}
 for S-fractions).
   \item[3)]
If $\bdelta$ is periodic of period 2
--- that is, $\delta_{2k-1} = x$ and $\delta_{2k} = y$ ---
then the sequence $(a_n)_{n \ge 0}$ generated by the T-fraction
is a linear transform of the sequence generated by the S-fraction
with the same coefficients $\balpha$:
see \cite[Propositions~3 and 15]{Barry_09} for some special cases,
and \cite{Sokal_totalpos} for the general case.
   \item[4)]
If $\delta_n = 0$ for all $n \ge 2$,
then the sequence $(a_n)_{n \ge 0}$ generated by the T-fraction
is a fairly simple nonlinear transform of the sequence generated by
the S-fraction with the same coefficients $\balpha$
(or alternatively the one with shifted coefficients $(\alpha_{n+1})_{n \ge 1}$):
this is a kind of ``renewal theory'' (see \cite{Sokal_totalpos}).
\end{itemize}
But T-fractions not of these four special types
are comparatively rare in the combinatorial literature.
It is our purpose here to provide a nontrivial example.

Two decades ago, Roblet and Viennot \cite{Roblet_96}
showed that the general T-fraction \reff{def.Ttype}
can be interpreted combinatorially as a generating function
for Dyck paths in which falls from peaks and non-peaks get different weights.
When $\bdelta = \bzero$ these two weights coincide,
and the Roblet--Viennot formula reduces to Flajolet's interpretation
of S-fractions in terms of Dyck paths.
More recently, several authors
\cite{Fusy_15,Oste_15,Josuat-Verges_18,Sokal_totalpos}
have independently found an alternate (and simpler)
combinatorial interpretation of the general T-fraction \reff{def.Ttype}:
namely, as a generating function for Schr\"oder paths with suitable
weights for each rise, fall and long level step.
When $\bdelta = \bzero$ the long level steps get zero weight,
and this formula again reduces to Flajolet's interpretation
of S-fractions in terms of Dyck paths.
We will review this Schr\"oder-path representation
in Section~\ref{subsec.prelim.1} below.

Our combinatorial example concerns the Ward polynomials $W_n(x)$,
which are defined as follows.
Firstly, the {\em Ward numbers}\/ $W(n,k)$
\cite{Ward_34} \cite[A134991/A181996/A269939]{OEIS}
are defined by the linear recurrence
\be
   W(n,k)  \;=\;  (n+k-1) \, W(n-1,k-1)  \:+\:  k \, W(n-1,k)
 \label{eq.ward.recurrence}
\ee
with initial condition $W(0,k) = \delta_{0k}$.
The triangular array of Ward numbers begins~as
%
%
\vspace*{-5mm}
\begin{table}[H]
\hspace*{-3mm}
\small
\begin{tabular}{c|rrrrrrrrr|r}
$n \setminus k$ & 0 & 1 & 2 & 3 & 4 & 5 & 6 & 7 & 8 & Row sums \\
\hline
0 & 1 &  &  &  &  &  &  &  &  & 1  \\
1 & 0 & 1 &  &  &  &  &  &  &  & 1  \\
2 & 0 & 1 & 3 &  &  &  &  &  &  & 4  \\
3 & 0 & 1 & 10 & 15 &  &  &  &  &  & 26  \\
4 & 0 & 1 & 25 & 105 & 105 &  &  &  &  & 236  \\
5 & 0 & 1 & 56 & 490 & 1260 & 945 &  &  &  & 2752  \\
6 & 0 & 1 & 119 & 1918 & 9450 & 17325 & 10395 &  &  & 39208  \\
7 & 0 & 1 & 246 & 6825 & 56980 & 190575 & 270270 & 135135 &  & 660032  \\
8 & 0 & 1 & 501 & 22935 & 302995 & 1636635 & 4099095 & 4729725 & 2027025 & 12818912  \\
\end{tabular}
\end{table}
\vspace*{-5mm}

\noindent
The row sums are \cite[A000311]{OEIS}.
It is easy to see that $W(n,n) = (2n-1)!!$
and that for $n \ge 1$ we have
$W(n,0) = 0$, $W(n,1) = 1$ and $W(n,2) = 2^{n+1} - n - 3$.
We then define the {\em Ward polynomials}\/
to be the row-generating polynomials of this triangular array:
\be
   W_n(x)  \;=\;  \sum_{k=0}^n W(n,k) \, x^k
   \;.
\ee
For some purposes it is convenient to use instead
the reversed Ward polynomials $\overline{W}_n(x) = x^n W_n(1/x)$,
which satisfy $\overline{W}_n(0) = (2n-1)!!$.
See \cite{Ward_34,Clark_99,Barbero_14,Barbero_15}
for further information on the Ward numbers and Ward polynomials.

Our first result is:

\begin{theorem}
   \label{thm1.1}
The ordinary generating function of the Ward polynomials
can be expressed by the T-fraction
\be
   \sum_{n=0}^\infty W_n(x) \, t^n
   \;=\;
   \cfrac{1}{1 - \cfrac{xt}{1 - t - \cfrac{2xt}{1 - 2t - \cfrac{3xt}{1- \cdots}}}}
  \label{eq.thm1.1}
\ee
with coefficients $\alpha_i = ix$ and $\delta_i = i-1$.
Equivalently, the ordinary generating function of the reversed Ward polynomials
can be expressed by the T-fraction
\be
   \sum_{n=0}^\infty \overline{W}_n(x) \, t^n
   \;=\;
   \cfrac{1}{1 - \cfrac{t}{1 - xt - \cfrac{2t}{1 - 2xt - \cfrac{3t}{1- \cdots}}}}
  \label{eq.thm1.1.reversed}
\ee
with coefficients $\alpha_i = i$ and $\delta_i = (i-1)x$.
\end{theorem}

\noindent
For the special case $x=1$, this T-fraction was known previously,
at least empirically \cite{Gladkovskii_13b}.

As preparation for the proof of Theorem~\ref{thm1.1},
let us mention two different combinatorial interpretations of the Ward numbers.
The first one, which is fairly well known,
is in terms of phylogenetic trees \cite{Steel_14}\footnote{
   Also called labelled hierarchies
   \cite{Leclerc_85,Flajolet_94} \cite[pp.~128--129 and 472--474]{Flajolet_09},
   Schr\"oder systems
   \cite{Schroder_1870,Comtet_70} \cite[pp.~223--224]{Comtet_74},
   Schr\"oder bracketings \cite{Riordan_76},
   or total partitions
   \cite[Example~5.2.5 and Exercise~5.43]{Stanley_99}.
},
the study of which goes back to Schr\"oder's fourth problem
\cite{Schroder_1870}.
A~{\em phylogenetic tree of type $(n,k)$}\/ is a rooted tree
that has $n+1$ labeled leaves (numbered $1,\ldots,n+1$)
and $k$ unlabeled internal vertices,
in which every internal vertex has at least two children.
Thus, a phylogenetic tree of type $(0,0)$ is just a single vertex labeled~1,
which is both the root and a leaf;
and it is easy to see that in all other cases we must have $1 \le k \le n$.
In fact, it is not difficult to see that
the number of phylogenetic trees of type $(n,k)$
is precisely the Ward number $W(n,k)$.\footnote{
   Note that the case $k=n$ corresponds to
   binary phylogenetic trees, i.e.\ those in which
   every internal vertex has exactly two children.
   The number of such trees is $W(n,n) = (2n-1)!!$
   \cite[Example~5.2.6]{Stanley_99}
   \cite[Section~2.5]{Callan_09}
   \cite{Diaconis_98,Prodinger_18}.
}
[{\sc Proof.}  The claim is clear for $n=0$ and $n=1$; now we use induction.
Leaf number $n+1$ can be connected either 
to an internal vertex on a phylogenetic tree of type $(n-1,k)$,
hence in $k W(n-1,k)$ ways;
or to an edge on a phylogenetic tree of type $(n-1,k-1)$,
creating a new internal vertex of out-degree~2,
hence in $(n+k-2) W(n-1,k-1)$ ways;
or to a new root, having the old root of
a phylogenetic tree of type $(n-1,k-1)$ as its other child,
hence in $W(n-1,k-1)$ ways.
This proves \reff{eq.ward.recurrence}.]
So the Ward polynomial $W_n(x)$
is the generating polynomial for phylogenetic trees on $n+1$ labeled leaves
in which each internal vertex gets a weight~$x$.

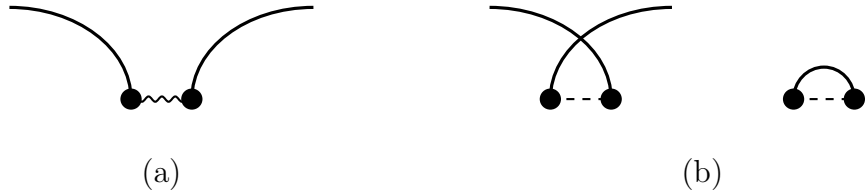
\begin{figure}[t]
\vspace*{1cm}
\begin{center}
  \begin{tikzpicture}[scale=0.8,line join=round]
    \PartArc{-1}{3}{0}{0.25}
    \PartArc{0}{4}{0.5}{0.25}
    \Wiggly{3}{4}
    \plotpermnobox[black] {}{0,0,1,1}
    \draw (3.5,0.2) node[below]{(a)};
  \end{tikzpicture}
  \hspace*{2cm}
  \begin{tikzpicture}[scale=0.8,line join=round]
    \PartArc{0}{4}{0}{0.25}
    \PartArc{-1}{3}{0.5}{0.25}
    \Dashed{3}{4}
    \Arc{7}{8}
    \Dashed{7}{8}
    \plotpermnobox[black] {}{0,0,1,1,0,0,1,1}
    \draw (5.5,0.2) node[below]{(b)};
  \end{tikzpicture}
\end{center}
\vspace*{-5mm}
\caption{(a)  A wiggly line in an augmented perfect matching.
         (b)  The two possibilities for a dashed line
              in a super-augmented perfect matching.}
\label{fig:augmented}
\end{figure}

The second interpretation, which seems to be new, is in terms
of what we shall call {\em augmented perfect matchings}\/.
Recall first that a perfect matching of $[2n] \eqdef \{1,\ldots,2n\}$
is simply a partition of $[2n]$ into $n$ pairs.
We call the smaller (resp.\ larger) element of each pair
an opener (resp.\ closer).
An augmented perfect matching is then a perfect matching
in which we may optionally also draw a wiggly line
on an edge $(i,i+1)$ whenever $i$ is a closer and $i+1$ is an opener
[see Figure~\ref{fig:augmented}(a)].
Then the number of augmented perfect matchings of $[2n]$
with $n-k$ wiggly lines
(or equivalently, $k$ closers that are not followed by a wiggly line)
is precisely the Ward number $W(n,k)$.
[{\sc Proof.} Let $M(n,\ell)$ be the the number of
  augmented perfect matchings of $[2n]$ with $\ell$ wiggly lines.
  The number of ways that the vertex $2n$ can be connected to
  a vertex $i \in [2n-1]$ that is not incident on a wiggly line
  is $(2n-1-\ell) M(n-1,\ell)$,
  because the vertex $i$ can be inserted in any of the gaps
  of an augmented perfect matching of $[2n-2]$ with $\ell$ wiggly lines,
  where ``gap'' means a space between adjacent vertices
  that are not joined by a wiggly line,
  or a space at the start or end.
  The number of ways that the vertex $2n$ can be connected to
  a vertex $i \in [2n-1]$ that is incident on a wiggly line
  is $(n-\ell) M(n-1,\ell-1)$,
  because the vertex $i$ together with its incident wiggly line
  (pointing to the left) can be inserted after any closer
  of an augmented perfect matching of $[2n-2]$ with $\ell-1$ wiggly lines
  that is not incident on a wiggly line.
  Comparing the recurrence
  $M(n,\ell) = (2n-1-\ell) M(n-1,\ell) + (n-\ell) M(n-1,\ell-1)$
  with \reff{eq.ward.recurrence},
  we see that $M(n,\ell) = W(n,n-\ell)$.]
So the reversed Ward polynomial $\overline{W}_n(x)$
is the generating polynomial for augmented perfect matchings of $[2n]$
in which each wiggly line gets a weight~$x$.
Equivalently, the Ward polynomial $W_n(x)$
is the generating polynomial for augmented perfect matchings of $[2n]$
in which each closer that is not incident on a wiggly line
gets a weight~$x$.

In Appendix~\ref{app.bijection} we construct a bijection between
augmented perfect matchings of $[2n]$ with $\ell$ wiggly lines
and phylogenetic trees of type $(n,n-\ell)$.

But we can go farther:
let us define a {\em super-augmented perfect matching}\/ of $[2n]$
to be a perfect matching of $[2n]$
in which we may optionally draw a wiggly line
on an edge $(i,i+1)$ whenever $i$ is a closer and $i+1$ is an opener,
and also optionally draw a dashed line 
on an edge $(i,i+1)$ whenever $i$ is an opener and $i+1$ is a closer;
however, it is not allowed for a wiggly line and a dashed line
to be incident on the same vertex.
Please note that for a pair $(i,i+1)$ connected by a dashed line,
there are two possibilities [see Figure~\ref{fig:augmented}(b)]:
the opener $i$ and the closer $i+1$ could belong to different arches
(which necessarily cross), or they could belong to the same arch.

Let us now say that a vertex is a {\em pure closer}\/
if it is a closer that is not incident on any wiggly or dashed line.
And let us say that a closer $k$ has {\em crossing number $m$}\/ ($m \ge 0$)
if there are $m$ quadruplets $i < j < k < l$
such that there are arches $(i,k)$ and $(j,l)$.

Now let $W_n(x,u,z,w',w'')$ be the
generating polynomial for super-augmented perfect matchings of $[2n]$
in which
each pure closer with crossing number 0 gets a weight~$x$,
each pure closer with crossing number $\ge 1$ gets a weight~$u$,
each dashed line for which the two endpoints belong to the same arch
gets a weight~$z$,
each other dashed line gets a weight~$w''$,
and each wiggly line gets a weight~$w'$.
Since each closer corresponds to precisely one of these five categories,
$W_n$ is a homogeneous polynomial of degree~$n$.
These five-variable polynomials have a T-fraction that generalizes
\reff{eq.thm1.1}:

\begin{theorem}
   \label{thm1.2}
The ordinary generating function of the polynomials $W_n(x,u,z,w',w'')$ 
can be expressed by the T-fraction
\be
\Scale[0.9]{
   \sum\limits_{n=0}^\infty W_n(x,u,z,w',w'') \, t^n
   \;=\;
   \cfrac{1}{1 - zt - \cfrac{xt}{1 - (z+w'+w'')t - \cfrac{(x+u)t}{1 - (z+2w'+2w'')t - \cfrac{(x+2u)t}{1- \cdots}}}}
}
  \label{eq.thm1.2}
\ee
with coefficients $\alpha_i = x + (i-1)u$ and $\delta_i = z + (i-1)(w'+w'')$.
\end{theorem}

\noindent
In particular, $W_n(x,u,z,w',w'')$ depends on $w'$ and $w''$
only via their sum $w \eqdef w' + w''$;
so we shall also write it as $W_n(x,u,z,w)$.
We call $W_n(x,u,z,w)$ the {\em generalized Ward polynomials}\/.
Note that $W_n(x) = W_n(x,x,0,1)$ and $\overline{W}_n(x) = W_n(1,1,0,x)$.
Let us also remark that $W_n(x,x,1,0)$ is \cite[A001498]{OEIS},
$W_n(x,x,1,1)$ is \cite[A112493]{OEIS},
and $W_n(1,1,z,1)$ is \cite[A298673]{OEIS}.

In Appendix~\ref{app.recurrence} we will show that
the sequence of four-variable polynomials
\linebreak
$W_n(x,u,z,w)$
satisfies a nonlinear differential recurrence,
and that the specialization to $x=u$
satisfies a linear differential recurrence.

But we can go much farther:
the main result of this paper is a ``master T-fraction''
that enumerates super-augmented perfect matchings
with an {\em infinite}\/ set of independent statistics
(Theorem~\ref{thm.matchings.Ttype.final1}).
This master T-fraction includes many combinatorially interesting polynomials
as special cases:
see Corollary~\ref{cor.matchings.18var} and the specializations
\reff{def.matching.12var}--\reff{def.weights.matchings.12var},
\reff{def.weights.matchings.12var.bis1}
and \reff{def.weights.matchings.12var.bis2}.
By further specialization of \reff{def.weights.matchings.12var.bis1},
we obtain Theorem~\ref{thm1.2}, which in turn implies Theorem~\ref{thm1.1}.

\bigskip

Let us conclude this introduction by mentioning briefly the connection
of the Ward numbers and Ward polynomials with other combinatorial
objects and problems:

\medskip

1)  The {\em 2-associated Stirling subset number}\/
$\assocstirlingsubset{n}{k}{2} \vphantom{\prod\limits_{}^\infty}$
is the number of partitions of an $n$-element set into $k$ blocks,
each of which has at least two elements;
by convention we set $\assocstirlingsubset{0}{k}{2} = \delta_{0k}$.
It is easy to see that these numbers satisfy the recurrence
\be
   \assocstirlingsubset{n}{k}{2}   \;=\;
      k \, \assocstirlingsubset{n-1}{k}{2}   \,+\,
      (n-1) \, \assocstirlingsubset{n-2}{k-1}{2}
   \qquad\hbox{for } n \ge 1
   \;.
 \label{eq.stirling.recurrence}
\ee
The first few $\assocstirlingsubset{n}{k}{2}$ are
\cite[pp.~221--222]{Comtet_74}
\cite[A008299/A137375]{OEIS}
%
%
\vspace*{-5mm}
\begin{table}[H]
\centering
\small
\begin{tabular}{c|rrrrrrrrrrr|r}
$n \setminus k$ & 0 & 1 & 2 & 3 & 4 & 5 & 6 & 7 & 8 & 9 & 10 & Row sums \\
\hline
0 & 1 &  &  &  &  &  &  &  &  &  &  & 1  \\
1 & 0 & 0 &  &  &  &  &  &  &  &  &  & 0  \\
2 & 0 & 1 & 0 &  &  &  &  &  &  &  &  & 1  \\
3 & 0 & 1 & 0 & 0 &  &  &  &  &  &  &  & 1  \\
4 & 0 & 1 & 3 & 0 & 0 &  &  &  &  &  &  & 4  \\
5 & 0 & 1 & 10 & 0 & 0 & 0 &  &  &  &  &  & 11  \\
6 & 0 & 1 & 25 & 15 & 0 & 0 & 0 &  &  &  &  & 41  \\
7 & 0 & 1 & 56 & 105 & 0 & 0 & 0 & 0 &  &  &  & 162  \\
8 & 0 & 1 & 119 & 490 & 105 & 0 & 0 & 0 & 0 &  &  & 715  \\
9 & 0 & 1 & 246 & 1918 & 1260 & 0 & 0 & 0 & 0 & 0 &  & 3425  \\
10 & 0 & 1 & 501 & 6825 & 9450 & 945 & 0 & 0 & 0 & 0 & 0 & 17722  \\
\end{tabular}
\end{table}
\vspace*{-5mm}
\noindent
Comparing \reff{eq.ward.recurrence} with \reff{eq.stirling.recurrence},
it is easy to see
that $W(n,k) = \assocstirlingsubset{n+k}{k}{2}$
\cite{Comtet_70} \cite[eqns.~(1.7) and (3.6)]{Carlitz_71}
\cite[pp.~6--7]{Riordan_76}.
This also has a nice bijective proof \cite{Comtet_70}:
given a phylogenetic tree of type $(n,k)$,
partition its $n+k$ non-root vertices into blocks consisting of
the sets of children of each of the $k$ internal vertices.\footnote{
   Here the leaves are labeled $1,\ldots,n+1$,
   but the internal vertices are {\em a priori}\/ unlabeled.
   Therefore, in order to make this construction precise,
   we must first use some algorithm to label the
   $k-1$ non-root internal vertices as $n+2,\ldots,n+k$
   (we can also label the root $n+k+1$ if we wish, but this plays no role).
   This labeling can be done in several different ways:
   see \cite{Comtet_70} \cite[Theorem~1]{Erdos_89}
   \cite[Lemma~5 and Theorem~5]{Haiman_89}
   \cite[Exercise~5.43, pp.~92 and 136--137]{Stanley_99}
   \cite[Theorem~2.1]{Gaiffi_15}.
   In each case one must prove that this indeed yields a bijection,
   i.e.\ that the phylogenetic tree can be reconstructed from
   the partition of $[n+k]$;  this is not completely trivial.
%
%
}
Note also that, under this bijection,
the number of internal vertices with $i$ children in the phylogenetic tree
equals the number of blocks of size $i$ in the partition,
as observed in
\cite[Theorem~1]{Erdos_89} \cite[Theorem~5]{Haiman_89}
\cite[Exercise~5.43, pp.~92 and 136--137]{Stanley_99}.


\medskip

2)  A {\em Stirling permutation}\/ \cite{Gessel_78} of order $n$ 
is a permutation $\sigma_1 \cdots \sigma_{2n}$
of the multiset $\{1,1,2,2,\ldots,n,n\}$
in which, for all $m \in [n]$, all numbers between the two occurrences
of $m$ are larger than $m$.
We say that an index $j \in [2n]$ is a {\em descent}\/
if $\sigma_j > \sigma_{j+1}$ or $j=2n$
(so by our definition the last index is always a descent).
The {\em second-order Eulerian number}\/ $\eulersecond{n}{k}$
is the number of Stirling permutations of order $n$
with exactly $k$ descents;
by convention we set $\eulersecond{0}{k} = \delta_{0k}$.\footnote{
   {\bf Warning:}  Our definition here of the second-order Eulerian numbers
   agrees with Gessel and Stanley \cite{Gessel_78}
   (who use the notation $B_{n,k}$)
   but disagrees with Graham, Knuth and Patashnik \cite[p.~270]{Graham_94}
   because we define the last index to be a descent while they do not:
   thus ${\eulersecond{0}{k}}^{\rm ours} = {\eulersecond{0}{k}}^{\rm GKP}
      = \delta_{0k}$
   but ${\eulersecond{n}{k}}^{\rm ours} = 
        \big\langle\!\! \euler{n}{k-1} \!\!\big\rangle^{\rm GKP}$
   for $n \ge 1$.
}
It is not difficult to see \cite[p.~27]{Gessel_78}
that these numbers satisfy the recurrence
\be
   \eulersecond{n}{k}   \;=\;
      (2n-k) \eulersecond{n-1}{k-1}  \,+\,
      k \eulersecond{n-1}{k}
   \qquad\hbox{for } n \ge 1
   \;.
 \label{eq.eulerian.recurrence}
\ee
The first few $\eulersecond{n}{k}$ are
\cite[A008517/A201637/A112007/A288874]{OEIS}
%
%
\vspace*{-5mm}
\begin{table}[H]
\centering
\small
\begin{tabular}{c|rrrrrrrrr|r}
$n \setminus k$ & 0 & 1 & 2 & 3 & 4 & 5 & 6 & 7 & 8 & Row sums \\
\hline
0 & 1 &  &  &  &  &  &  &  &  & 1  \\
1 & 0 & 1 &  &  &  &  &  &  &  & 1  \\
2 & 0 & 1 & 2 &  &  &  &  &  &  & 3  \\
3 & 0 & 1 & 8 & 6 &  &  &  &  &  & 15  \\
4 & 0 & 1 & 22 & 58 & 24 &  &  &  &  & 105  \\
5 & 0 & 1 & 52 & 328 & 444 & 120 &  &  &  & 945  \\
6 & 0 & 1 & 114 & 1452 & 4400 & 3708 & 720 &  &  & 10395  \\
7 & 0 & 1 & 240 & 5610 & 32120 & 58140 & 33984 & 5040 &  & 135135  \\
8 & 0 & 1 & 494 & 19950 & 195800 & 644020 & 785304 & 341136 & 40320 & 2027025  \\
\end{tabular}
\end{table}
\vspace*{-5mm}
\noindent
It is easy to see from \reff{eq.eulerian.recurrence}
that $\eulersecond{n}{n} = n!$
and $\sum_{k=0}^n \eulersecond{n}{k} = (2n-1)!!$.
Define now the {\em second-order Eulerian polynomials}\/
\be
   E_n^{[2]}(x)
   \;=\;
   \sum_{k=0}^n \eulersecond{n}{k} \, x^k
\ee
and the reversed polynomials
$\overline{E}_n^{[2]}(x) = x^n E_n^{[2]}(1/x)$;
note that $\overline{E}_n^{[2]}(0) = n!$
and $\overline{E}_n^{[2]}(1) = (2n-1)!!$.
By manipulating the recurrences \reff{eq.ward.recurrence}
and \reff{eq.eulerian.recurrence},
it is not difficult to show that
\be
   \overline{W}_n(x)  \;=\;  \overline{E}_n^{[2]}(1+x)
 \label{eq.ward.euler2}
\ee
or equivalently
\begin{subeqnarray}
   E_n^{[2]}(x)  & = &  (1-x)^n \, W_n \Bigl( {x \over 1-x} \Bigr)  \\[2mm]
   W_n(x)  & = &  (1+x)^n \, E_n^{[2]} \Bigl( {x \over 1+x} \Bigr)
\end{subeqnarray}
In particular we have $\overline{W}_n(-1) = n!$
\cite[Theorem~3.1 and Section~4]{Callan_15}.

The identity \reff{eq.ward.euler2} also has a nice combinatorial
interpretation and proof.  In a perfect matching of $[2n]$,
let us say that a pair $(i,i+1)$ is a closer/opener pair
if $i$~is a closer and $i+1$ is an opener.
Note that counting augmented perfect matchings
with a weight~$x$ for each wiggly line
--- as $\overline{W}_n(x)$ does ---
is equivalent to counting perfect matchings
with a weight~$1+x$ for each closer/opener pair.
Now let $M'(n,\ell)$ be the number of perfect matchings of $[2n]$
with $\ell$ closer/opener pairs.
We claim that $M'(n,\ell) = \eulersecond{n}{n-\ell}$.
[{\sc Proof.}
Write $\clop(\pi)$ for the number of closer/opener pairs
in a perfect matching $\pi$.
Now consider a perfect matching $\pi$ of $[2n]$,
and suppose that the vertex $2n$ is paired with $i \in [2n-1]$;
let $\pi'$ be the perfect matching of $[2n-2]$ obtained from $\pi$
by deleting the arch $(i,2n)$ and renumbering vertices.
If $i=1$ or $i-1$ is an opener, then $\clop(\pi') = \clop(\pi)$.
If $i \ge 2$ and $i-1$ is a closer, then $(i-1,i)$ is a closer/opener pair;
if also $i+1$ is an opener, then $\clop(\pi') = \clop(\pi)$,
otherwise $\clop(\pi') = \clop(\pi) - 1$.
If $\clop(\pi) = \ell$,
these insertions can be done in $n$, $\ell$ and
$n-\ell$ ways, respectively.
So $M'(n,\ell) = (n+\ell) M'(n-1,\ell) + (n-\ell) M'(n-1,\ell-1)$,
which is equivalent to \reff{eq.eulerian.recurrence}.]\footnote{
   Note also that $M'(n,0) = n!$ has an easy proof:
   if there are no closer/opener pairs,
   then the elements $1,\ldots,n$ are openers
   and $n+1,\ldots,2n$ are closers,
   and a perfect matching is obtained by pairing these elements
   in one of the $n!$ possible ways.
}
So the reversed second-order Eulerian polynomial $\overline{E}_n^{[2]}(x)$
is the generating polynomial for perfect matchings of $[2n]$
in which each closer/opener pair gets a weight~$x$.
This interpretation seems to be new.


\medskip

3)  There is also a very interesting multivariate generalization
of the Ward polynomials.
As explained above, the Ward polynomial $W_n(x)$
is the generating polynomial for phylogenetic trees on $n+1$ labeled leaves
in which each internal vertex gets a weight~$x$.
Now let $\bfx = (x_i)_{i \ge 1}$ be an infinite collection of indeterminates,
and let $\bW_n(\bfx) = \bW_n(x_1,\ldots,x_n)$
be the generating polynomial for phylogenetic trees on $n+1$ labeled leaves
in which each internal vertex with $i$ ($\ge 2$) children
gets a weight~$x_{i-1}$.
We call this the {\em multivariate Ward polynomial}\/;
it is quasi-homogeneous of degree~$n$ when $x_i$ is given weight~$i$.
Thus $W_n(x) = \bW_n(x,\ldots,x)$;
and from the quasi-homogeneity one deduces
$\overline{W}_n(x) = \bW_n(1,x,x^2,\ldots,x^{n-1})$.
The first few $\bW_n$ are \cite[A134685]{OEIS}
\begin{subeqnarray}
   \bW_0  & = &  1  \\
   \bW_1  & = &  x_1  \\
   \bW_2  & = &  3 x_1^2 + x_2  \\
   \bW_3  & = &  15 x_1^3 + 10 x_1 x_2 + x_3  \\
   \bW_4  & = &  105 x_1^4 + 105 x_1^2 x_2 + 15 x_1 x_3 + 10 x_2^2 + x_4 
 \label{eq.multiward.0-4}
\end{subeqnarray}
A straightforward argument
(generalizing \cite[pp.~128--129]{Flajolet_09}
 or \cite[Examples~5.2.5 and 5.2.6]{Stanley_99})
shows that the exponential generating function
\be
   \scrw(t;\bfx)
   \;\eqdef\;
   \sum_{n=0}^\infty \bW_n(\bfx) \, {t^{n+1} \over (n+1)!}
   \;=\;
   t \,+\, \sum_{n=2}^\infty \bW_{n-1}(\bfx) \, {t^n \over n!}
 \label{def.scrw}
\ee
satisfies the functional equation
\be
   \scrw(t;\bfx)
   \;=\;
   t \,+\, \sum_{n=2}^\infty x_{n-1} \, {\scrw(t;\bfx)^n \over n!}
   \;,
 \label{eq.functeqn.scrw}
\ee
where the term $n$ in the sum corresponds to the case
in which the root has $n$ children.
In other words, $\scrw(t;\bfx)$ is the compositional inverse of
the generic power series
\be
   F(t;\bfx) \;\eqdef\; t \,-\, \sum_{n=2}^\infty x_{n-1} \, {t^n \over n!}
   \;,
 \label{def.Ftx}
\ee
which therefore satisfies
[substituting $t \leftarrow F(t;\bfx)$ in \reff{def.scrw}]
\be
   F(t;\bfx)
   \;=\;
   t \,-\, \sum_{n=2}^\infty \bW_{n-1}(\bfx) \, {F(t;\bfx)^n \over n!}
   \;.
 \label{eq.functeqn.Ftx}
\ee
In this context of series inversion,
the polynomials $\bW_n(x_1,\ldots,x_n)$ can be found in the books of
Riordan \cite[p.~181]{Riordan_68}\footnote{
   With the erratum given by Riordan in \cite[p.~7]{Riordan_76}.
   Riordan \cite[p.~7]{Riordan_76} also observes
   (attributing this remark to Neil Sloane)
   that the specialization of $\bW_n(x_1,\ldots,x_n)$
   [defined via compositional inverses] to $x_1 = \ldots = x_n = x$
   enumerates phylogenetic trees by number of internal vertices.
   But his formula $S_n(y) = Z_n(y,\ldots,y)$ [his $Z_n$ is our $\bW_n$]
   should read $S_n(y) = Z_{n-1}(y,\ldots,y)$,
   since his $S_n(y)$ enumerates phylogenetic trees on $n$ leaves, not $n+1$.
}
and Comtet \cite[p.~151]{Comtet_74},
but without the interpretation in~terms of phylogenetic trees.
See also \cite[A134685]{OEIS} \cite{Haiman_89} \cite[Theorem~2]{Drake_07}
for further information concerning the multivariate Ward polynomials.

Since $\bW_n(x_1,\ldots,x_n) = x_n \,+\,$
a polynomial with integer coefficients in $x_1,\ldots,x_{n-1}$,
any sequence $\ba = (a_n)_{n \ge 0}$ with $a_0 = 1$ in any commutative ring $R$
can be written as a specialization of $\bW_n(x_1,\ldots,x_n)$
with suitable coefficients $(x_i)_{i \ge 1}$ in $R$.
For instance, we have
\begin{subeqnarray}
   x_1  & = &  a_1  \\
   x_2  & = &  -3 a_1^2 + a_2  \\
   x_3  & = &  15 a_1^3 + 10 a_1 a_2 + a_3  \\
   x_4  & = &  -105 a_1^4 + 105 a_1^2 a_2 - 15 a_1 a_3 - 10 a_2^2 + a_4 
 \label{eq.multiward.inverse.0-4}
\end{subeqnarray}
It will be observed that this is identical to \reff{eq.multiward.0-4}
with some sign changes (namely, a minus sign on each monomial
with an even number of factors $a_i$):
that is, $-x_n = \bW_n(-a_1,\ldots,-a_n)$.
This is in fact easy to prove: if we write
\be
   \scra(t)
   \;\eqdef\;
   \sum_{n=0}^\infty a_n \, {t^{n+1} \over (n+1)!}
   \;=\;
   t \,+\, \sum_{n=2}^\infty a_{n-1} \, {t^n \over n!}
   \;,
 \label{def.scra}
\ee
then setting $a_n = \bW_n(x_1,\ldots,x_n)$ in \reff{eq.functeqn.scrw}
yields the functional equation
\be
   \scra(t)
   \;=\;
   t \,+\, \sum_{n=2}^\infty x_{n-1} \, {\scra(t)^n \over n!}
   \;.
 \label{eq.functeqn.scra}
\ee
Comparing \reff{def.scra}/\reff{eq.functeqn.scra}
with \reff{def.Ftx}/\reff{eq.functeqn.Ftx}
and recalling that $\bW_{n-1}(\bfx) = a_{n-1}$,
we see that they are related by $(a_n,x_n) \to (-x_n,-a_n)$,
which proves the claim.

In particular, taking $a_n$ equal to the generalized Ward polynomial
$W_n(x,u,z,w)$, we have
\begin{subeqnarray}
   x_1  & = &  x + z  \\
   x_2  & = &  ux + wx - x^2 - 3xz - 2z^2 \\
   x_3  & = &  3 u^2 x + 4 u w x - 3 u x^2 - 5 u x z + w^2 x - 4 w x^2
     - 6 w x z + 5 x^2 z + 11 x z^2 + 6 z^3
     \nonumber \\[-1mm] \\[-6mm]
        & \vdots &  \nonumber
 \label{eq.multiward.inverse.Wnxuzw}
\end{subeqnarray}
We have not yet succeeded in guessing the general formula for these $x_i$;
but in the special case $u=x$
we can show (see Appendix~\ref{app.recurrence}) that
\be
   x_{n-1}
   \;=\;
   (-1)^n \, (n-1)! \, \Bigl( 1 + {x \over w} \Bigr) \, z^{n-1}
     \:+\:
   {x \over w} \, \prod_{j=1}^{n-1} (w - jz)
 \label{eq.multiward.inverse.Wnxuzw.u=x}
\ee
(note that this is a polynomial in $x,z,w$ since the terms in $1/w$ cancel),
corresponding to
\be
   F(t;\bfx)
   \;=\;
   {1 + {x \over w} \over z} \, \log(1 + zt)
     \:-\:
   {x \over w^2} \, \bigl[ (1 + zt)^{w/z} \,-\, 1 \bigr]
   \;.
 \label{eq.multiward.inverse.Wnxuzw.u=x.F}
\ee
When $z=0$ and $w=1$ this reduces to $x_i = x$
and $F(t;\bfx) = t - x(e^t - 1 - t)$.
It is curious that the $x_i$ in \reff{eq.multiward.inverse.Wnxuzw}
and \reff{eq.multiward.inverse.Wnxuzw.u=x}
are {\em not}\/ all nonnegative (even pointwise in $x,u,z,w \ge 0$),
even though the output polynomials
$\bW_n(x_1,\ldots,x_n) = W_n(x,u,z,w)$
are of course coefficientwise nonnegative in $x,u,z,w$.

It would be interesting to extend our work on the Ward polynomials
to this multivariate generalization.
%
%

\medskip

4) Finally, let us comment on the implications of our results
for the theory of Hankel-total positivity
\cite{Sokal_flajolet,Sokal_totalpos},
which was in fact our original motivation for studying these T-fractions.
If $\bP = (P_n(\bfx))_{n \ge 0}$
is a sequence of polynomials with real coefficients
in one or more indeterminates $\bfx$,
we say that this sequence is {\em coefficientwise Hankel-totally positive}\/
if every minor of the infinite Hankel matrix
$H_\infty(\bP) = (P_{i+j}(\bfx))_{i,j \ge 0}$
is a polynomial with nonnegative coefficients.
One fundamental result of this theory is \cite{Sokal_totalpos}:
whenever the ordinary generating function
$\sum_{n=0}^\infty P_n(\bfx) \, t^n$
is given by a T-fraction \reff{def.Ttype}
where all the coefficients $\alpha_i$ and $\delta_i$
are polynomials in $\bfx$ with nonnegative coefficients,
the sequence $\bP = (P_n(\bfx))_{n \ge 0}$ is
coefficientwise Hankel-totally positive.
It follows from this that all the T-fractions derived in this paper
--- including the most general one, Theorem~\ref{thm.matchings.Ttype.final1}
--- are coefficientwise Hankel-totally positive in the appropriate variables.
In particular, the Ward polynomials $W_n(x)$
and the generalized Ward polynomials $W_n(x,u,z,w)$
are coefficientwise Hankel-totally positive.

But now a mystery arises from the relation \reff{eq.ward.euler2}
that connects the (reversed) Ward polynomials
with the (reversed) second-order Eulerian polynomials.
It follows from \reff{eq.ward.euler2} and \reff{eq.thm1.1.reversed}
that the ordinary generating function
of the reversed second-order Eulerian polynomials
is given by a T-fraction with coefficients
$\delta_i = (i-1)(x-1)$ and $\alpha_i = i$.\footnote{
   This T-fraction is also a consequence of
   \cite[Section~5.2, item~(7)]{Callan_09}
   together with the Roblet--Viennot \cite{Roblet_96}
   interpretation of T-fractions as a generating function
   for Dyck paths with special weights for peaks.
}
Here $\delta_i$ is {\em not}\/ coefficientwise nonnegative
(or even pointwise nonnegative when $0 \le x < 1$),
so the general theory \cite{Sokal_totalpos}
says nothing about the Hankel-total positivity of
the reversed second-order Eulerian polynomials.
And yet, we find empirically that the sequence of 
reversed second-order Eulerian polynomials
{\em is}\/ coefficientwise Hankel-totally positive:
we have tested this through the $13 \times 13$ Hankel matrix.
This total positivity (if indeed it is true) is, alas, completely unexplained.
By \reff{eq.ward.euler2} it {\em implies}\/
the coefficientwise Hankel-total positivity
of the reversed Ward polynomials --- which we know to be true ---
but not conversely.
In our opinion, proving the coefficientwise Hankel-total positivity
of the (reversed) second-order Eulerian polynomials
is the major problem left open by our work.

{\bf Note added:}
This mystery has now been resolved \cite{latpath_SRTR}.
The reversed second-order Eulerian polynomials $\overline{E}_n^{[2]}(x)$
are {\em also}\/ given by a 2-branched S-fraction \cite{latpath_SRTR}
with coefficients
$\balpha = (\alpha_i)_{i \ge 2} = 1,1,x,2,2,2x,3,3,3x,\ldots\:$.
Since these $\alpha_i$ are coefficientwise nonnegative,
the general theory of branched S-fractions \cite{latpath_SRTR}
implies the coefficientwise Hankel-total positivity.

\bigskip

The plan of this paper is as follows:
In Section~\ref{sec.statement} we state our results
on the enumeration of super-augmented perfect matchings:
these include a very general T-fraction with infinitely many indeterminates
(Theorem~\ref{thm.matchings.Ttype.final1})
as well as numerous special cases that count statistics of combinatorial
interest.
In Section~\ref{sec.prelim} we recall some basic facts
concerning the combinatorial interpretation of continued fractions
and the concept of labeled Schr\"oder paths.
In Section~\ref{sec.proofs} we prove
Theorem~\ref{thm.matchings.Ttype.final1}
by constructing a bijection from the set
of super-augmented perfect matchings
onto a set of labeled 2-colored Schr\"oder paths.
In Appendix~\ref{app.bijection} we construct a bijection between
augmented perfect matchings and phylogenetic trees.
In Appendix~\ref{app.recurrence} we show that
the generalized Ward polynomials $W_n(x,u,z,w)$
satisfy a nonlinear differential recurrence.

\section{Statement of results}   \label{sec.statement}

We begin by reviewing some S-fractions for perfect matchings
that were recently obtained by Zeng and one of us
\cite{Sokal-Zeng_masterpoly}.
We then state our generalizations,
which are T-fractions for super-augmented perfect matchings.
Let us stress that none of our results in this paper
depend on the results of \cite{Sokal-Zeng_masterpoly};
we give here self-contained proofs of all our generalizations.

\subsection{S-fractions for perfect matchings}  \label{subsec.statement.S}

Euler showed \cite[section~29]{Euler_1760} that the
generating function of the odd semifactorials can be represented
as an S-type continued fraction
\be
   \sum_{n=0}^\infty (2n-1)!! \: t^n
   \;=\;
   \cfrac{1}{1 - \cfrac{1t}{1 - \cfrac{2t}{1 - \cfrac{3t}{1- \cdots}}}}
 \label{eq.2n-1semifact.contfrac}
\ee
with coefficients $\alpha_n = n$.\footnote{
   See also \cite[Section~2.6]{Callan_09} for a combinatorial proof
   of \reff{eq.2n-1semifact.contfrac}
   based on a counting of height-labeled Dyck paths.
}
Since $(2n-1)!!$ enumerates perfect matchings of a $2n$-element set,
it is natural to seek polynomial refinements of this sequence
that enumerate perfect matchings of $[2n]$ according to some
natural statistic(s).
Note that we can regard a perfect matching either as a special type
of set partition (namely, one in which all blocks are of size 2)
or as a special type of permutation (namely, a fixed-point-free involution).
We will adopt here the former interpretation,
and write $\pi \in \scrm_{2n} \subseteq \Pi_{2n}$.
If $i,j \in [2n]$ are paired in the perfect matching $\pi$,
we write $i \sim_\pi j$
(or just $i \sim j$ if $\pi$ is clear from the context).


Inspired by \reff{eq.2n-1semifact.contfrac},
let us introduce \cite{Sokal-Zeng_masterpoly}
the polynomials $M_n(x,y,u,v)$ defined by the continued fraction
\be
   \sum_{n=0}^\infty M_n(x,y,u,v) \: t^n
   \;=\;
   \cfrac{1}{1 - \cfrac{xt}{1 - \cfrac{(y+v)t}{1 - \cfrac{(x+2u)t}{1- \cfrac{(y+3v)t}{1 - \cfrac{(x+4u)t}{1 - \cfrac{(y+5v)t}{1-\cdots}}}}}}}
 \label{eq.matching.fourvar.contfrac}
\ee
with coefficients
\begin{subeqnarray}
   \alpha_{2k-1}  & = &  x + (2k-2) u \\
   \alpha_{2k}    & = &  y + (2k-1) v
 \label{def.weights.matching.fourvar}
\end{subeqnarray}
Clearly $M_n(x,y,u,v)$ is a homogeneous polynomial of degree~$n$.
Since $M_n(1,1,1,1) = (2n-1)!!$,
it is natural to expect that $M_n(x,y,u,v)$ enumerates
perfect matchings of $[2n]$ according to some natural trivariate statistic.
Let us now explain, following \cite{Sokal-Zeng_masterpoly},
what this statistic is.

Let $\pi$ be a perfect matching of $[2n]$.
We recall that a vertex $i \in [2n]$ is called an opener (resp.\ closer)
if it is the smaller (resp.\ larger) element of its pair.
Let us say that an opener $j$ (with corresponding closer $k$)
is a {\em record}\/ if there does not exist an opener $i < j$
that is paired with a closer $l > k$.
In other words, $j$ is a record if there does not exist
an arch $(i,l)$ that nests above the arch $(j,k)$.
Similarly, let us say that a closer $k$ (with corresponding opener $j$)
is an {\em antirecord}\/ if there does not exist a closer $l > k$
that is paired with an opener $i < j$.
In other words, $k$ is an antirecord if and only if $j$ is a record.\footnote{
   This terminology of records and antirecords comes from the
   alternate interpretation of perfect matchings as fixed-point-free
   involutions.
   In general, for a permutation $\sigma \in \Sym_n$,
   an index $i \in [n]$ is called a
   {\em record}\/ (or {\em left-to-right maximum}\/)
   if $\sigma(j) < \sigma(i)$ for all $j < i$,
   and is called an {\em antirecord}\/ (or {\em right-to-left minimum}\/)
   if $\sigma(j) > \sigma(i)$ for all $j > i$.
   See \cite{Sokal-Zeng_masterpoly} for further discussion
   of record-antirecord statistics in permutations.
}

Let us now classify the closers of $\pi$ into four types:
\begin{itemize}
   \item {\em even closer antirecord}\/ (ecar)\ \ 
      [i.e.\ $i$ is even and is an antirecord];
   \item {\em odd closer antirecord}\/ (ocar)\ \ 
      [i.e.\ $i$ is odd and is an antirecord];
   \item {\em even closer non-antirecord}\/ (ecnar)\ \ 
      [i.e.\ $i$ is even and is not an antirecord];
   \item {\em odd closer non-antirecord}\/ (ocnar)\ \ 
      [i.e.\ $i$ is odd and is not an antirecord].
\end{itemize}
Similarly, we classify the openers of $\pi$ into four types:
\begin{itemize}
   \item {\em even opener record}\/ (eor)\ \ 
      [i.e.\ $i$ is even and is a record];
   \item {\em odd opener record}\/ (oor)\ \ 
      [i.e.\ $i$ is odd and is a record];
   \item {\em even opener non-record}\/ (eonr)\ \ 
      [i.e.\ $i$ is even and is not a record];
   \item {\em odd opener non-record}\/ (oonr)\ \ 
      [i.e.\ $i$ is odd and is not a record].
\end{itemize}
Then Sokal and Zeng \cite{Sokal-Zeng_masterpoly} proved that
the polynomials $M_n(x,y,u,v)$ defined by
\reff{eq.matching.fourvar.contfrac}/\reff{def.weights.matching.fourvar}
have the combinatorial interpretation
\begin{subeqnarray}
   M_n(x,y,u,v)
   & = &
   \sum_{\pi \in \scrm_{2n}}
      x^{\ecar(\pi)} y^{\ocar(\pi)}
         u^{\ecnar(\pi)}  v^{\ocnar(\pi)}
                        \slabel{eq.matching.fourvar.a}  \\[2mm]
   & = &
   \sum_{\pi \in \scrm_{2n}}
      x^{\oor(\pi)} y^{\eor(\pi)}
         u^{\oonr(\pi)}  v^{\eonr(\pi)}
                        \slabel{eq.matching.fourvar.b}
   \;.
 \label{eq.matching.fourvar}
\end{subeqnarray}
The interpretations \reff{eq.matching.fourvar.a} and
\reff{eq.matching.fourvar.b} are of course trivially equivalent
under the involution $i \to 2n+1-i$,
which preserves the structure of a perfect matching
but interchanges even with odd, opener with closer, and record with antirecord.

Sokal and Zeng \cite{Sokal-Zeng_masterpoly} also generalized
the four-variable polynomial $M_n(x,y,u,v)$
by adding weights for crossings and nestings.
%
Given a perfect matching $\pi \in \scrm_{2n}$,
we say that a quadruplet  $i < j < k < l$ forms a
{\em crossing}\/ if $i \sim_\pi k$ and $j \sim_\pi l$,
and a {\em nesting}\/ if $i \sim_\pi l$ and $j \sim_\pi k$.
We write $\crr(\pi)$ and $\nee(\pi)$ for the numbers of
crossings and nestings of $\pi$.
We now define the six-variable polynomial
\begin{subeqnarray}
   M_n(x,y,u,v,p,q)
   & = &
   \sum_{\pi \in \scrm_{2n}}
      x^{\ecar(\pi)} y^{\ocar(\pi)}
         u^{\ecnar(\pi)}  v^{\ocnar(\pi)}
         p^{\crr(\pi)} q^{\nee(\pi)}
     \qquad
      \\[2mm]
   & = &
   \sum_{\pi \in \scrm_{2n}}
      x^{\oor(\pi)} y^{\eor(\pi)}
         u^{\oonr(\pi)}  v^{\eonr(\pi)}
         p^{\crr(\pi)} q^{\nee(\pi)}
   \;.
 \label{def.matching.fourvar.pq}
\end{subeqnarray}
Here the equality of (\ref{def.matching.fourvar.pq}a)
and (\ref{def.matching.fourvar.pq}b) follows again from
the involution $i \to 2n+1-i$,
which preserves the numbers of crossings and nestings.
Sokal and Zeng \cite{Sokal-Zeng_masterpoly} showed that
the ordinary generating function of the polynomials $M_n(x,y,u,v,p,q)$
has an S-fraction that generalizes
\reff{eq.matching.fourvar.contfrac}/\reff{def.weights.matching.fourvar},
namely
\be
   \sum_{n=0}^\infty M_n(x,y,u,v,p,q) \: t^n
   \;\,=\;\,
   \cfrac{1}{1 - \cfrac{xt}{1 - \cfrac{(py+qv)t}{1 - \cfrac{(p^2 x+q \, [2]_{p,q} u)t}{1- \cfrac{(p^3 y+ q \, [3]_{p,q} v)t}{1 - \cdots}}}}}
   \label{eq.thm.matching.fourvar.pq.Stype}
\ee
with coefficients
\begin{subeqnarray}
   \alpha_{2k-1}  & = &  p^{2k-2} x + q \, [2k-2]_{p,q} u \\
   \alpha_{2k}    & = &  p^{2k-1} y + q \, [2k-1]_{p,q} v
 \label{def.weights.matching.fourvar.pq.Stype}
\end{subeqnarray}
where
\be
   [n]_{p,q}
   \;\eqdef\;
   {p^n - q^n \over p-q}
   \;=\;
   \sum\limits_{j=0}^{n-1} p^j q^{n-1-j}
 \label{def.npq}
\ee
for an integer $n \ge 0$.
When $p=q=1$ this reduces to
\reff{eq.matching.fourvar.contfrac}/\reff{def.weights.matching.fourvar}.
Note that if $u=x$ and/or $v=y$,
then the weights \reff{def.weights.matching.fourvar.pq.Stype} simplify to
$\alpha_{2k-1} = [2k-1]_{p,q} \, x$ and $\alpha_{2k} = [2k]_{p,q} \, y$,
respectively.
For the special case $x=y=u=v$,
the S-fraction \reff{eq.thm.matching.fourvar.pq.Stype}
was obtained previously by Kasraoui and Zeng \cite{Kasraoui_06}
(see also \cite[p.~3280]{Blitvic_12}).

But Sokal and Zeng \cite{Sokal-Zeng_masterpoly} went even farther,
by defining statistics that count the number of quadruplets
$i < j < k < l$ that form crossings or nestings
{\em with a particular vertex $k$ in third position}\/:
\begin{eqnarray}
      \crr(k,\pi)
   & = &
   \#\{ i<j<k<l \colon\: i \sim_\pi k \hbox{ and } j \sim_\pi l \}
       \label{def.cr.k.pi}
         \\[2mm]
   \nee(k,\pi)
   & = &
   \#\{ i<j<k<l \colon\: i \sim_\pi l \hbox{ and } j \sim_\pi k \}
       \label{def.ne.k.pi}
\end{eqnarray}
Note that $\crr(k,\pi)$ and $\nee(k,\pi)$ can be nonzero
only when $k$ is a closer (and $k \ge 3$);
in geometrical terms, $\crr(k,\pi)$ [resp.\ $\nee(k,\pi)$]
is the number of arches of $\pi$ that cross (resp.\ nest over)
the arch whose right endpoint is $k$.
We obviously have
\begin{subeqnarray}
   \crr(\pi)
   & = &
   \!\!
   \sum\limits_{k \in {\rm closers}}  \!\! \crr(k,\pi)
         \\[2mm]
   \nee(\pi)
   & = &
   \!\!
   \sum\limits_{k \in {\rm closers}}  \!\! \nee(k,\pi)
 \label{eq.crne.total}
\end{subeqnarray}

Sokal and Zeng \cite{Sokal-Zeng_masterpoly} also defined
\be
   \qne(k,\pi)
   \;=\;
   \#\{ i<k<l \colon\:  i \sim_\pi l  \}
 \label{def.qne}
\ee
(they call this a {\em quasi-nesting}\/ of the vertex $k$).
We will use this statistic only when $k$ is an opener:
in this case $\qne(k,\pi)$ counts the number of times that the opener $k$
occurs in {\em second}\/ position in a crossing or nesting:
when $(i,l) \in \scrg_\pi$ is a pair contributing to $\qne(k,\pi)$,
the quadruplet ${i<k<l,m}$ [where $(k,m) \in \scrg_\pi$]
must be either a crossing or nesting (according as $m>l$ or $m<l$),
but we do not keep track of which one it is.

Now introduce two infinite families of indeterminates
$\bsfa = (\sfa_\ell)_{\ell \ge 0}$ and
$\bsfb = (\sfb_{\ell,\ell'})_{\ell,\ell' \ge 0}$,
and define the polynomials $M_n(\bsfa,\bsfb)$ by
\be
   M_n(\bsfa,\bsfb)
   \;=\;
   \sum_{\pi \in \scrm_{2n}}
   \;
   \prod\limits_{i \in {\rm openers}}  \!\! \sfa_{\qne(i,\pi)}
   \prod\limits_{i \in {\rm closers}}  \! \sfb_{\crr(i,\pi),\, \nee(i,\pi)}
   \;.
 \label{def.Mn.ab}
\ee
Sokal and Zeng \cite{Sokal-Zeng_masterpoly} showed that
the ordinary generating function of the polynomials
$M_n(\bsfa,\bsfb)$ has the S-fraction
\be
   \sum_{n=0}^\infty M_n(\bsfa,\bsfb) \: t^n
   \;=\;
   \cfrac{1}{1 -  \cfrac{\sfa_{0} \sfb_{00} t}{1 -  \cfrac{\sfa_1 (\sfb_{01} + \sfb_{10}) t}{1 - \cfrac{\sfa_2 (\sfb_{02} + \sfb_{11} + \sfb_{20}) t}{1 - \cdots}}}}
   \label{eq.thm.matchings.Stype.final1}
\ee
with coefficients
\be
   \alpha_n   \;=\;   \sfa_{n-1} \, \sfb^\star_{n-1}
 \label{def.weights.matchings.Stype.final1}
\ee
where
\be
   \sfb^\star_{n-1}  \;\eqdef\;  \sum_{\ell=0}^{n-1} \sfb_{\ell,n-1-\ell}
   \;.
 \label{def.bstar0}
\ee
This ``master S-fraction for perfect matchings''
contains
\reff{eq.matching.fourvar.contfrac}--\reff{eq.matching.fourvar}
and
\reff{def.matching.fourvar.pq}--\reff{def.weights.matching.fourvar.pq.Stype}
as special cases \cite{Sokal-Zeng_masterpoly}.

\medskip

\begin{quote}
{\bf Warning:}
The definitions \reff{def.cr.k.pi}/\reff{def.ne.k.pi} and \reff{def.Mn.ab}
and the S-fraction \reff{eq.thm.matchings.Stype.final1}--\reff{def.bstar0}
are taken from an early draft of \cite{Sokal-Zeng_masterpoly}.
The final version of \cite{Sokal-Zeng_masterpoly} applies a reversal
$i \mapsto 2n+1-i$ to these definitions and results:
that is, crossings and nestings are defined with a particular vertex $j$
in {\em second}\/ rather than third position;
the two-index indeterminate in \reff{def.Mn.ab}
is associated with openers rather than closers;
and the letters $\sfa$ and $\sfb$ are interchanged.
But the versions stated here are of course correct.
\end{quote}

\subsection{T-fractions for super-augmented perfect matchings}
    \label{subsec.statement.T}

Recall that a {\em super-augmented perfect matching}\/ of $[2n]$
is a perfect matching of $[2n]$
in which we may optionally draw a wiggly line
on an edge $(i,i+1)$ whenever $i$ is a closer and $i+1$ is an opener,
and also optionally draw a dashed line 
on an edge $(i,i+1)$ whenever $i$ is an opener and $i+1$ is a closer;
however, it is not allowed for a wiggly line and a dashed line
to be incident on the same vertex.
We denote a super-augmented perfect matching of $[2n]$
by $\tau \in \scrm^\star_{2n}$,
and we write $\pi(\tau) \in \scrm_{2n}$ for the underlying perfect matching.
Of course, we say that $i$ is an opener (resp.\ closer) in $\tau$
in case it is an opener (resp.\ closer) in $\pi(\tau)$.

For $\tau \in \scrm^\star_{2n}$,
we call a vertex $i$ a
\begin{itemize}
   \item {\em pure opener}\/ if it is an opener
      and not incident on any wiggly or dashed line;
   \item {\em pure closer}\/ if it is a closer
      and not incident on any wiggly or dashed line;
   \item {\em wiggly closer}\/ if it is a closer
      that is incident on a wiggly line;
   \item {\em dashed closer}\/ if it is a closer
      that is incident on a dashed line.
\end{itemize}
We also call a pair $(i,i+1)$ a
\begin{itemize}
   \item {\em wiggly pair}\/ if $i$ is a closer, $i+1$ is an opener,
      and they are connected by a wiggly line;
   \item {\em dashed pair}\/ if $i$ is an opener, $i+1$ is a closer,
      and they are connected by a dashed line.
\end{itemize}

The statistics \reff{def.cr.k.pi}--\reff{def.qne}
apply equally well to super-augmented perfect matchings
if we evaluate them on $\pi = \pi(\tau)$.
We will generalize \reff{def.Mn.ab},
not by introducing any new statistics,
but simply by refining the distinctions among types of vertices.
We introduce four infinite families of indeterminates
$\bsfa = (\sfa_\ell)_{\ell \ge 0}$,
$\bsfb = (\sfb_{\ell,\ell'})_{\ell,\ell' \ge 0}$,
$\bsff = (\sff_{\ell,\ell'})_{\ell,\ell' \ge 0}$ and
$\bsfg = (\sfg_{\ell,\ell'})_{\ell,\ell' \ge 0}$,
and define the polynomials $M_n(\bsfa,\bsfb,\bsff,\bsfg)$ by
\begin{eqnarray}
   M_n(\bsfa,\bsfb,\bsff,\bsfg)
   & = & 
   \sum_{\pi \in \scrm^\star_{2n}}
   \;
   \prod\limits_{i \in {\rm pureop}}  \!\!\! \sfa_{\qne(i,\pi)}
   \prod\limits_{i \in {\rm purecl}}  \!\! \sfb_{\crr(i,\pi),\, \nee(i,\pi)}
            \nonumber \\[2mm]
   &   &
   \hspace{8.7mm}
   \prod\limits_{(i,i+1) \in {\rm wiggly}}  \!\!\!\!\! \sff_{\crr(i,\pi),\, \nee(i,\pi)}
   \prod\limits_{(i,i+1) \in {\rm dashed}}  \!\!\!\!\! \sfg_{\crr(i+1,\pi),\, \nee(i+1,\pi)}
   \;, \qquad
 \label{def.Mn}
\end{eqnarray}
where pureop (resp.\ purecl) denotes the pure openers (resp.\ pure closers).
Our main result is then the following:

\begin{theorem}[Master T-fraction for super-augmented perfect matchings]
   \label{thm.matchings.Ttype.final1}
The ordinary generating function of the polynomials
$M_n(\bsfa,\bsfb,\bsff,\bsfg)$ has the T-type continued fraction
\be
   \sum_{n=0}^\infty M_n(\bsfa,\bsfb,\bsff,\bsfg) \: t^n
   \;=\;
   \cfrac{1}{1 - \sfg_{00} t - \cfrac{\sfa_{0} \sfb_{00} t}{1 - (\sff_{00} + \sfg_{01} + \sfg_{10})t - \cfrac{\sfa_1 (\sfb_{01} + \sfb_{10}) t}{1 - \cdots}}}
   \label{eq.thm.matchings.Ttype.final1}
\ee
with coefficients
\begin{subeqnarray}
   \alpha_n  & = &   \sfa_{n-1} \, \sfb^\star_{n-1}
      \slabel{def.weights.matchings.Ttype.final1} \\[2mm]
   \delta_n  & = &   \sff^\star_{n-2} \,+\, \sfg^\star_{n-1}
   \label{eq.thm.matchings.Ttype.final1.coeffs}
\end{subeqnarray}
where
\be
 \sfb^\star_{n-1}  \;\eqdef\;  \sum_{\ell=0}^{n-1} \sfb_{\ell,n-1-\ell}
   \label{def.bstar}
\ee
and likewise for $\sff$ and $\sfg$.
\end{theorem}

\noindent
When $\sff = \sfg  = 0$, wiggly and dashed lines are forbidden,
and Theorem~\ref{thm.matchings.Ttype.final1}
reduces to the S-fraction
\reff{eq.thm.matchings.Stype.final1}/\reff{def.weights.matchings.Stype.final1}.
We will prove Theorem~\ref{thm.matchings.Ttype.final1}
in Section~\ref{sec.proofs}.

Let us now show some specializations of
Theorem~\ref{thm.matchings.Ttype.final1}
that generalize the S-fractions shown in the preceding subsection.
We begin by defining a 18-variable polynomial
$M_n(x,y,u,v,x',y',u',v',x'',y'',u'',v'',p,q,p',q',p'',q'')$
that generalizes (\ref{def.matching.fourvar.pq}a)
by distinguishing the three types of closers: pure, wiggly and dashed.
We let $\ecar^\circ(\tau)$, $\ecar'(\tau)$ and $\ecar''(\tau)$
be the number of even closer antirecords in each of these three classes,
and likewise for $\ocar$, $\ecnar$ and $\ocnar$.
We also distinguish three types of crossings and nestings,
according to the nature of the closer that occurs in third position:
we thus let $\crr^\circ(\tau)$, $\crr'(\tau)$ and $\crr''(\tau)$
be the sum of $\crr(k,\pi)$ over vertices $k$ that are,
respectively, pure closers, wiggly closers and dashed closers.
We then define:
\begin{eqnarray}
   & & \!\!\!\!\!\!
   M_n(x,y,u,v,x',y',u',v',x'',y'',u'',v'',p,q,p',q',p'',q'')
          \nonumber \\[2mm]
   & & \qquad
   \;=\;
   \sum_{\tau \in \scrm^\star_{2n}}
      \:
      x^{\ecar^\circ(\tau)} y^{\ocar^\circ(\tau)}
         u^{\ecnar^\circ(\tau)}  v^{\ocnar^\circ(\tau)}
             \nonumber \\[-2mm]
   &   &  \:\hspace*{2.6cm}
      (x')^{\ecar'(\tau)} (y')^{\ocar'(\tau)}
         (u')^{\ecnar'(\tau)}  (v')^{\ocnar'(\tau)}
             \nonumber \\[2mm]
   &   &  \:\hspace*{2.6cm}
      (x'')^{\ecar''(\tau)} (y'')^{\ocar''(\tau)}
         (u'')^{\ecnar''(\tau)}  (v'')^{\ocnar''(\tau)}
             \nonumber \\[2mm]
   &   &  \:\hspace*{2.6cm}
         p^{\crr^\circ(\tau)} q^{\nee^\circ(\tau)}
         (p')^{\crr'(\tau)} (q')^{\nee'(\tau)}
         (p'')^{\crr''(\tau)} (q'')^{\nee''(\tau)}
   \qquad\qquad
 \label{def.matching.18var}
\end{eqnarray}
We have:

\begin{corollary}
   \label{cor.matchings.18var}
The ordinary generating function of the polynomials
\reff{def.matching.18var}
has the T-type continued fraction
\begin{eqnarray}
   & & \hspace*{-9mm}
   \sum_{n=0}^\infty 
      M_n(x,y,u,v,x',y',u',v',x'',y'',u'',v'',p,q,p',q',p'',q'') \: t^n
         \nonumber \\[2mm]
   & &
   \!\!\!\!\!\!=\;
\Scale[0.75]{
   \cfrac{1}{1 - x'' t - \cfrac{xt}{1 - (x' + p'' y'' + q'' v'')t - \cfrac{(py+qv)t}{1 - (p' y' + q' v' + (p'')^2 x'' + q'' \, [2]_{p,1} u'')t - \cfrac{(p^2 x+q \, [2]_{p,q} u)t}{1- \cdots}}}}
}
  \qquad
   \label{eq.cor.matchings.18var}
\end{eqnarray}
with coefficients
\begin{subeqnarray}
   \alpha_{2k-1}  & = &  p^{2k-2} x + q \, [2k-2]_{p,q} u \\[1mm]
   \alpha_{2k}    & = &  p^{2k-1} y + q \, [2k-1]_{p,q} v \\[1mm]
   \delta_1       & = &  x''  \\[1mm]
   \delta_{2k-1}  & = &  (p')^{2k-3} y' + q' \, [2k-3]_{p',q'} v' + (p'')^{2k-2} x'' + q'' \, [2k-2]_{p'',q''} u'' \; \hbox{ for } k \ge 2
             \nonumber \\[-1mm] \\
   \delta_{2k}    & = &  (p')^{2k-2} x' + q' \, [2k-2]_{p',q'} u' + (p'')^{2k-1} y'' + q'' \, [2k-1]_{p'',q''} v''
 \label{def.weights.matchings.18var}
\end{subeqnarray}
where $[n]_{p,q}$ is defined in \reff{def.npq}.
\end{corollary}

When $x'=y'=u'=v'=x''=y''=u''=v''=0$
(thereby forbidding wiggly and dashed lines),
Corollary \ref{cor.matchings.18var} reduces to the S-fraction
\reff{eq.thm.matching.fourvar.pq.Stype}/%
\reff{def.weights.matching.fourvar.pq.Stype}.

The proof of Corollary \ref{cor.matchings.18var} will be based
on the following easy combinatorial lemma \cite{Sokal-Zeng_masterpoly}:

\begin{lemma}[Closers in perfect matchings]
   \label{lemma.matchings.pi}
Let $\pi \in \scrm_{2n}$ be a perfect matching of $[2n]$,
and let $k \in [2n]$ be a closer of $\pi$.  Then:
\begin{itemize}
  \item[(a)]  $k$ has the same parity as $\crr(k,\pi) + \nee(k,\pi)$.
  \item[(b)]  $k$ is an antirecord if and only if $\nee(k,\pi) = 0$.
\end{itemize}
\end{lemma}

\proof
(a) Let $j$ be the opener that is paired with the closer $k$.
For each $i < k$, let $\sigma(i)$ ($\neq i$)
be the element with which it is paired.
Then the set $\{i \colon\; i < k \}$, which has cardinality $k-1$,
can be partitioned as
\be
   \{j\}
   \;\cup\;
   \{i < k \colon\: \sigma(i) < k \}
   \;\cup\;
   \{i < j \colon\: \sigma(i) > k \}
   \;\cup\;
   \{j < i < k \colon\: \sigma(i) > k \}
   \;.
\ee
The first of these sets has cardinality 1;
the second has even cardinality;
the third has cardinality $\nee(k,\pi)$;
and the fourth has cardinality $\crr(k,\pi)$.

(b) As explained earlier, a closer $k$ (with corresponding opener $j$)
is an antirecord if and only if there does not exist a closer $l > k$
that is paired with an opener $i < j$.
But this is precisely the statement that $\nee(k,\pi) = 0$.
\qed

\proofof{Corollary~\ref{cor.matchings.18var}}
In \reff{def.Mn} we set $\sfa_\ell = 1$ for all $\ell \ge 0$
and
\begin{subeqnarray}
   \sfb_{\ell,\ell'}
   & = &
   p^\ell q^{\ell'} \times
   \begin{cases}
       x     &  \textrm{if $\ell' = 0$ and $\ell$ is even}  \\
       y     &  \textrm{if $\ell' = 0$ and $\ell$ is odd}   \\
       u     &  \textrm{if $\ell' \ge 1$ and $\ell+\ell'$ is even}
                                                            \\
       v     &  \textrm{if $\ell' \ge 1$ and $\ell+\ell'$ is odd}
   \end{cases}
      \\[4mm]
   \sff_{\ell,\ell'}
   & = &
   (p')^\ell \, (q')^{\ell'} \times
   \begin{cases}
       x'    &  \textrm{if $\ell' = 0$ and $\ell$ is even}  \\
       y'    &  \textrm{if $\ell' = 0$ and $\ell$ is odd}   \\
       u'    &  \textrm{if $\ell' \ge 1$ and $\ell+\ell'$ is even}
                                                            \\
       v'    &  \textrm{if $\ell' \ge 1$ and $\ell+\ell'$ is odd}
   \end{cases}
      \\[4mm]
   \sfg_{\ell,\ell'}
   & = &
   (p'')^\ell \, (q'')^{\ell'} \times
   \begin{cases}
       x''    &  \textrm{if $\ell' = 0$ and $\ell$ is even}  \\
       y''    &  \textrm{if $\ell' = 0$ and $\ell$ is odd}   \\
       u''    &  \textrm{if $\ell' \ge 1$ and $\ell+\ell'$ is even}
                                                            \\
       v''    &  \textrm{if $\ell' \ge 1$ and $\ell+\ell'$ is odd}
   \end{cases}
 \label{eq.defs.bfg}
\end{subeqnarray}
By Lemma~\ref{lemma.matchings.pi}(a,b) and
the definitions of $\crr^\circ(\pi)$, etc.,
we obtain the polynomial \reff{def.matching.18var}.
Then
\be
   \sfb^\star_{n-1}  \;\eqdef\;  \sum_{\ell=0}^{n-1} \sfb_{\ell,n-1-\ell}
                     \;=\;
   \begin{cases}
       p^{n-1} x + q [n-1]_{p,q} u     &  \textrm{if $n$ is odd}  \\[1mm]
       p^{n-1} y + q [n-1]_{p,q} v     &  \textrm{if $n$ is even}
   \end{cases}
\ee
and similarly for $\sff$ and $\sfg$.
The result then follows from Theorem~\ref{thm.matchings.Ttype.final1}.
\qed

If we further specialize \reff{def.matching.18var}
by setting $x=y$, $u=v$, $x'=y'$, $u'=v'$, $x''=y''$, $u''=v''$,
we obtain 12-variable polynomials that count closers of each of three types
(pure, wiggly or dashed)
according to whether they are antirecords or not
--- but now forgetting whether they are even or odd ---
and also count crossings and nestings of each of three types
(pure, wiggly or dashed):
\begin{eqnarray}
   & & \!\!\!\!\!\!
   M_n(x,x,u,u,x',x',u',u',x'',x'',u'',u'',p,q,p',q',p'',q'')
          \nonumber \\[2mm]
   & &
   \;=\;
   \sum_{\tau \in \scrm^\star_{2n}}
      \:
      x^{\car^\circ(\tau)} u^{\cnar^\circ(\tau)}
      (x')^{\car'(\tau)} (u')^{\cnar'(\tau)}
      (x'')^{\car''(\tau)} (u'')^{\cnar''(\tau)}
             \nonumber \\[-2mm]
   &   &  \:\hspace*{1.85cm}
         p^{\crr^\circ(\tau)} q^{\nee^\circ(\tau)}
         (p')^{\crr'(\tau)} (q')^{\nee'(\tau)}
         (p'')^{\crr''(\tau)} (q'')^{\nee''(\tau)}
   \qquad\qquad
 \label{def.matching.12var}
\end{eqnarray}
Their T-fraction is
\begin{eqnarray}
   & & \hspace*{-9mm}
   \sum_{n=0}^\infty 
      M_n(x,x,u,u,x',x',u',u',x'',x'',u'',u'',p,q,p',q',p'',q'') \: t^n
         \nonumber \\[2mm]
   & &
   \!\!\!\!\!\!=\;
\Scale[0.72]{
   \cfrac{1}{1 - x'' t - \cfrac{xt}{1 - (x' + p'' x'' + q'' u'')t - \cfrac{(px+qu)t}{1 - (p' x' + q' u' + (p'')^2 x'' + q'' \, [2]_{p'',q''} u'')t - \cfrac{(p^2 x+q \, [2]_{p,q} u)t}{1- \cdots}}}}
}
  \qquad
   \label{eq.cor.matchings.12var}
\end{eqnarray}
with coefficients
\begin{subeqnarray}
   \alpha_n  & = &  p^{n-1} x + q \, [n-1]_{p,q} u \\[1mm]
   \delta_1  & = &  x''  \\[1mm]
   \delta_n  & = &  (p')^{n-2} x' + q' \, [n-2]_{p',q'} u' + (p'')^{n-1} x'' + q'' \, [n-1]_{p'',q''} u'' \; \hbox{ for } n \ge 2 \qquad\qquad
 \label{def.weights.matchings.12var}
\end{subeqnarray}
Note that the forms of the coefficients \reff{def.weights.matchings.12var}
are no longer alternating between even and odd,
because we are no longer distinguishing even closers from odd closers.

If we further specialize this to $u' = x'$,
then the coefficients \reff{def.weights.matchings.12var} simplify to
\begin{subeqnarray}
   \alpha_n  & = &  p^{n-1} x + q \, [n-1]_{p,q} u \\[1mm]
   \delta_n  & = &  [n-1]_{p',q'} x' + (p'')^{n-1} x'' + q'' \, [n-1]_{p'',q''} u''
 \label{def.weights.matchings.12var.bis1}
\end{subeqnarray}
And if we specialize \reff{def.weights.matchings.12var}
or \reff{def.weights.matchings.12var.bis1}
to $u=x$ and/or $u''=x''$, then the combinations involving those variables
simplify as well:
\begin{subeqnarray}
   p^{n-1} x + q \, [n-1]_{p,q} u   & \to &    [n]_{p,q} x   \qquad \\[2mm]
   (p'')^{n-1} x'' + q'' \, [n-1]_{p'',q''} u''
                                    & \to &    [n]_{p'',q''} x'' \qquad
 \label{def.weights.matchings.12var.bis2}
\end{subeqnarray}
When the three specializations $u=x$, $u'=x'$, $u''=x''$ are made,
we are no longer distinguishing antirecords from non-antirecords.

It is now an easy matter to deduce Theorem~\ref{thm1.2}.
We recall that the polynomial $W_n(x,u,z,w',w'')$ was defined
in the Introduction to be the generating polynomial
for super-augmented perfect matchings of $[2n]$ in which
each pure closer with crossing number 0 gets a weight~$x$,
each pure closer with crossing number $\ge 1$ gets a weight~$u$,
each dashed line for which the two endpoints belong to the same arch
gets a weight~$z$,
each other dashed line gets a weight~$w''$,
and each wiggly line gets a weight~$w'$.
We now observe that whenever $i$ is an opener and $i+1$ is a closer,
we have $\crr(i+1,\pi) = 0$ if and only if there is an arch from $i$ to $i+1$:
for if there is an arch from $i$ to $i+1$,
then obviously it is not crossed by any other arch;
and if there is not an arch from $i$ to $i+1$,
then there must exist $j < i < i+1 < k$
such that there are arches $(j,i+1)$ and $(i,k)$.
This applies in particular whenever there is a dashed line from $i$ to $i+1$.
So saying that
each dashed line for which the two endpoints belong to the same arch
gets a weight~$z$
and each other dashed line gets a weight~$w''$
is equivalent to saying that each dashed closer
with crossing number 0 (resp.~$\ge 1$) gets a weight~$z$ (resp.~$w''$).
In other words, we are evaluating \reff{def.Mn} at
\begin{subeqnarray}
   \sfa_\ell  & = &  1  \\[2mm]
   \sfb_{\ell,\ell'}
   & = &
   \begin{cases}
       x     &  \textrm{if $\ell = 0$}  \\
       u     &  \textrm{if $\ell \ge 1$}
   \end{cases}
      \\[2mm]
   \sff_{\ell,\ell'}   & = &   w'   \\[2mm]
   \sfg_{\ell,\ell'}
   & = &
   \begin{cases}
       z    &  \textrm{if $\ell = 0$}  \\
       w''  &  \textrm{if $\ell \ge 1$}
   \end{cases}
 \label{eq.defs.abfg.thm1.2}
\end{subeqnarray}
Then \reff{eq.thm.matchings.Ttype.final1.coeffs} becomes
\begin{subeqnarray}
   \alpha_n  & = &  x + (n-1) u \\[1mm]
   \delta_n  & = &  z + (n-1)(w'+w'')
\end{subeqnarray}
as claimed in Theorem~\ref{thm1.2}.

\section{Preliminaries for the proofs}   \label{sec.prelim}

Our proof of Theorem~\ref{thm.matchings.Ttype.final1}
will be based on Flajolet's \cite{Flajolet_80}
combinatorial interpretation of continued fractions
in terms of Dyck and Motzkin paths,
adapted slightly to handle Schr\"oder paths,
together with a bijection that maps
super-augmented perfect matchings to labeled Schr\"oder paths.
We begin by reviewing briefly these two ingredients.

\subsection{Combinatorial interpretation of continued fractions}
   \label{subsec.prelim.1}

Recall that a {\em Motzkin path}\/ of length $n \ge 0$
is a path $\omega = (\omega_0,\ldots,\omega_n)$
in the right quadrant $\N \times \N$,
starting at $\omega_0 = (0,0)$ and ending at $\omega_n = (n,0)$,
whose steps $s_j = \omega_j - \omega_{j-1}$
are $(1,1)$ [``rise''], $(1,-1)$ [``fall''] or $(1,0)$ [``level''].
We write $\scrm_n$ for the set of Motzkin paths of length~$n$,
and $\scrm = \bigcup_{n=0}^\infty \scrm_n$.
A Motzkin path is called a {\em Dyck path}\/ if it has no level steps.
A Dyck path always has even length;
we write $\scrd_{2n}$ for the set of Dyck paths of length~$2n$,
and $\scrd = \bigcup_{n=0}^\infty \scrd_{2n}$.

Let ${\bf a} = (a_i)_{i \ge 0}$, ${\bf b} = (b_i)_{i \ge 1}$
and ${\bf c} = (c_i)_{i \ge 0}$ be indeterminates;
we will work in the ring $\Z[[{\bf a},{\bf b},{\bf c}]]$
of formal power series in these indeterminates.
To each Motzkin path $\omega$ we assign a weight
$W(\omega) \in \Z[{\bf a},{\bf b},{\bf c}]$
that is the product of the weights for the individual steps,
where a rise starting at height~$i$ gets weight~$a_i$,
a~fall starting at height~$i$ gets weight~$b_i$,
and a level step at height~$i$ gets weight~$c_i$.
Flajolet \cite{Flajolet_80} showed that
the generating function of Motzkin paths
can be expressed as a continued fraction:

\begin{theorem}[Flajolet's master theorem]
   \label{thm.flajolet}
We have
\be
   \sum_{\omega \in \scrm}  W(\omega)
   \;=\;
   \cfrac{1}{1 - c_0 - \cfrac{a_0 b_1}{1 - c_1 - \cfrac{a_1 b_2}{1- c_2 - \cfrac{a_2 b_3}{1- \cdots}}}}
 \label{eq.thm.flajolet}
\ee
as an identity in $\Z[[{\bf a},{\bf b},{\bf c}]]$.
\end{theorem}

In particular, if $a_{i-1} b_i = \beta_i t^2$ and $c_i = \gamma_i t$
(note that the parameter $t$ is conjugate to the length of the Motzkin path),
we have
\be
   \sum_{n=0}^\infty t^n \sum_{\omega \in \scrm_n}  W(\omega)
   \;=\;
   \cfrac{1}{1 - \gamma_0 t - \cfrac{\beta_1 t^2}{1 - \gamma_1 t - \cfrac{\beta_2 t^2}{1 - \cdots}}}
   \;\,,
 \label{eq.flajolet.motzkin}
\ee
so that the generating function of Motzkin paths with height-dependent weights
is given by the J-type continued fraction \reff{def.Jtype}.
Similarly, if $a_{i-1} b_i = \alpha_i t$ and $c_i = 0$
(note that $t$ is now conjugate to the semi-length of the Dyck path), we have
\be
   \sum_{n=0}^\infty t^n \sum_{\omega \in \scrd_{2n}}  W(\omega)
   \;=\;
   \cfrac{1}{1 - \cfrac{\alpha_1 t}{1 - \cfrac{\alpha_2 t}{1 - \cdots}}}
   \;\,,
 \label{eq.flajolet.dyck}
\ee
so that the generating function of Dyck paths with height-dependent weights
is given by the S-type continued fraction \reff{def.Stype}.

Let us now show how to handle Schr\"oder paths within this framework.
A {\em Schr\"oder path}\/ of length $2n$ ($n \ge 0$)
is a path $\omega = (\omega_0,\ldots,\omega_{2n})$
in the right quadrant $\N \times \N$,
starting at $\omega_0 = (0,0)$ and ending at $\omega_{2n} = (2n,0)$,
whose steps are $(1,1)$ [``rise''], $(1,-1)$ [``fall'']
or $(2,0)$ [``long level''].
We write $s_j$ for the step starting at abscissa $j-1$.
If the step $s_j$ is a rise or a fall,
we set $s_j = \omega_j - \omega_{j-1}$ as before.
If the step $s_j$ is a long level step,
we set $s_j = \omega_{j+1} - \omega_{j-1}$ and leave $\omega_j$ undefined;
furthermore, in this case there is no step $s_{j+1}$.
We write $h_j$ for the height of the Schr\"oder path at abscissa $j$
whenever this is defined, i.e.\ $\omega_j = (j,h_j)$.
Please note that $\omega_{2n} = (2n,0)$ and $h_{2n} = 0$
are always well-defined,
because there cannot be a long level step starting at abscissa $2n-1$.
We write $\scrs_{2n}$ for the set of Schr\"oder paths of length~$2n$,
and $\scrs = \bigcup_{n=0}^\infty \scrs_{2n}$.

There is an obvious bijection between Schr\"oder paths and Motzkin paths:
namely, every long level step is mapped onto a level step.
If we apply Flajolet's master theorem with
$a_{i-1} b_i = \alpha_i t$ and $c_i = \delta_{i+1} t$
to the resulting Motzkin path
(note that $t$ is now conjugate to the semi-length
 of the underlying Schr\"oder path),
we obtain
\be
   \sum_{n=0}^\infty t^n \sum_{\omega \in \scrs_{2n}}  W(\omega)
   \;=\;
   \cfrac{1}{1 - \delta_1 t - \cfrac{\alpha_1 t}{1 - \delta_2 t - \cfrac{\alpha_2 t}{1 - \cdots}}}
   \;\,,
 \label{eq.flajolet.schroder}
\ee
so that the generating function of Schr\"oder paths
with height-dependent weights
is given by the T-type continued fraction \reff{def.Ttype}.
More precisely, every rise gets a weight~1,
every fall starting at height~$i$ gets a weight $\alpha_i$,
and every long level step at height~$i$ gets a weight $\delta_{i+1}$.
This combinatorial interpretation of T-fractions in terms of Schr\"oder paths
was found recently by several authors
\cite{Fusy_15,Oste_15,Josuat-Verges_18,Sokal_totalpos}.

\subsection{Labelled Schr\"oder paths}

Many authors, starting with Flajolet \cite{Flajolet_80},
have used bijections from combinatorial objects
onto labeled Motzkin or Dyck paths
in order to prove J-fraction or S-fraction expansions
for the (weighted) ordinary generating functions of those objects.
Here we will do the same with labeled Schr\"oder paths
in order to prove T-fraction expansions.
The definitions are as follows:

Let ${\bf A} = (A_k)_{k \ge 0}$, ${\bf B} = (B_k)_{k \ge 1}$
and ${\bf C} = (C_k)_{k \ge 0}$ be sequences of nonnegative integers.
An {\em $({\bf A},{\bf B},{\bf C})$-labeled Schr\"oder path of length $2n$}\/
is a pair $(\omega,\xi)$
where $\omega = (\omega_0,\ldots,\omega_{2n})$
is a Schr\"oder path of length $2n$,
and $\xi = (\xi_1,\ldots,\xi_{2n})$ is a sequence of integers satisfying
\be
   1  \:\le\: \xi_i  \:\le\:
   \begin{cases}
       A(h_{i-1})  & \textrm{if step $s_i$ is a rise (starting at height $h_{i-1}$)}
              \\[1mm]
       B(h_{i-1})  & \textrm{if step $s_i$ is a fall (starting at height $h_{i-1}$)}
              \\[1mm]
       C(h_{i-1})  & \textrm{if step $s_i$ is a long level step (at height $h_{i-1}$)}
   \end{cases}
 \label{eq.xi.ineq}
\ee
[For typographical clarity
 we have here written $A(k)$ as a synonym for $A_k$, etc.]
If step $s_i$ is undefined (because step $s_{i-1}$ was a long level step),
then $\xi_i$ is also undefined.
We denote by $\scrs_{2n}({\bf A},{\bf B},{\bf C})$
the set of $({\bf A},{\bf B},{\bf C})$-labeled Schr\"oder paths
of length $2n$.

Let us stress that the numbers $A_k$, $B_k$ and $C_k$ are allowed
to take the value 0.
Whenever this happens, the path $\omega$ is forbidden to take a step
of the specified kind at the specified height.

We shall also make use of multicolored Schr\"oder paths.
A {\em $k$-colored Schr\"oder path}\/ is simply a Schr\"oder path
in which each long level step has been given a ``color''
from the set $\{1,2,\ldots,k\}$.
In other words, we distinguish $k$ different types of long level steps.
An {\em $({\bf A},{\bf B},{\bf C}^{(1)},\ldots,{\bf C}^{(k)})$-labeled
 $k$-colored Schr\"oder path of length $2n$}\/
is then defined in the obvious way,
where we use the sequence ${\bf C}^{(j)}$ to bound
the label $\xi_i$ when step $i$ is a long level step of type $j$.
We denote by $\scrs_{2n}({\bf A},{\bf B},{\bf C}^{(1)},\ldots,{\bf C}^{(k)})$
the set of $({\bf A},{\bf B},{\bf C}^{(1)},\ldots,{\bf C}^{(k)})$-labeled
$k$-colored Schr\"oder paths of length $2n$.

\section{Proof of Theorem~\ref{thm.matchings.Ttype.final1}} \label{sec.proofs}

We will prove Theorem~\ref{thm.matchings.Ttype.final1}
by constructing a bijection from the set $\scrm_{2n}^\star$
of super-augmented perfect matchings of $[2n]$
onto the set of $({\bf A},{\bf B},{\bf C}^{(1)},{\bf C}^{(2)})$-labeled
2-colored Schr\"oder paths of length $2n$, where
\begin{subeqnarray}
   A_k        & = &  1         \quad\qquad\hbox{for $k \ge 0$}  \\
   B_k        & = &  k         \quad\qquad\hbox{for $k \ge 1$}  \\
   C_k^{(1)}  & = &  k         \quad\qquad\hbox{for $k \ge 0$}  \\
   C_k^{(2)}  & = &  k+1       \quad\:\hbox{for $k \ge 0$}
 \label{def.abc.KZ}
\end{subeqnarray}
When restricted to ordinary perfect matchings of $[2n]$
(i.e.\ super-augmented perfect matchings with no wiggly or dashed lines),
our bijection maps onto labeled Dyck paths
and coincides with the bijection used by Flajolet \cite{Flajolet_80}
and Kasraoui--Zeng \cite{Kasraoui_06}.\footnote{
   More precisely,
   Flajolet \cite{Flajolet_80} and Kasraoui--Zeng \cite{Kasraoui_06}
   defined bijections of set partitions of $[n]$
   onto labeled Motzkin paths of length $n$.
   When restricted to perfect matchings of $[2n]$
   (i.e.\ set partitions of $[2n]$ in which every block has cardinality 2),
   their bijections map onto labeled Dyck paths of length $2n$.
   The Flajolet and Kasraoui--Zeng bijections for set partitions
   are slightly different, but they coincide when restricted
   to perfect matchings.
}

We will begin by explaining how the Schr\"oder path $\omega$ is defined;
then we will explain how the labels $\xi$ are defined;
next we will prove that the mapping is indeed a bijection;
next we will translate the various statistics from
$\scrm_{2n}^\star$ to our labeled Schr\"oder paths;
and finally we will sum over labels $\xi$ to obtain the weight $W(\omega)$
associated to a Schr\"oder path $\omega$,
which upon applying \reff{eq.flajolet.schroder}
will yield Theorem~\ref{thm.matchings.Ttype.final1}.

\bigskip

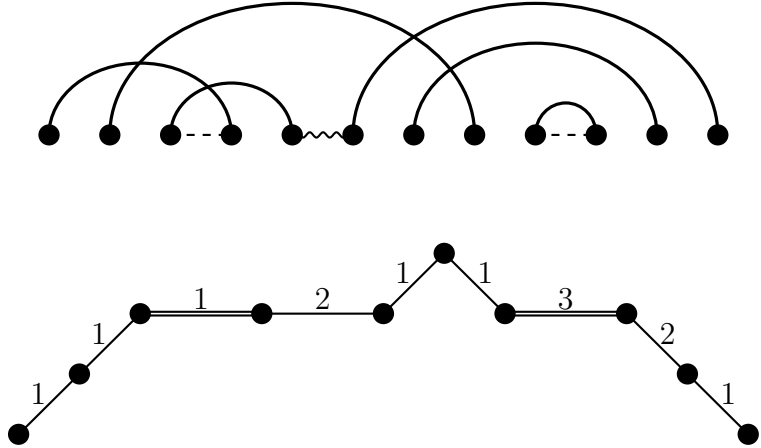
\begin{figure}[t]
\vspace*{1cm}
\begin{center}
  \begin{tikzpicture}[scale=0.8,line join=round]
     \Arc{1}{4}
     \Arc{2}{8}
     \Arc{3}{5}
     \Arc{6}{12}
     \Arc{7}{11}
     \Arc{9}{10}
    \Wiggly{5}{6}
    \Dashed{3}{4}
    \Dashed{9}{10}
    
    \plotpermnobox[black] {}{1,1,1,1,1,1,1,1,1,1,1,1}
  \end{tikzpicture}
\end{center}
\begin{center}
\begin{tikzpicture}[scale=0.8,line join=round]
\node at (2,5) {$~$};
\draw [black,thick] (1,1)--(2,2);
    \draw [black,thick] (2,2)--(3,3);
    \draw [black,thick,double] (3,3)--(5,3);
    \draw [black,thick] (5,3)--(7,3);
    \draw [black,thick] (7,3)--(8,4);
    \draw [black,thick] (8,4)--(9,3);
    \draw [black,thick,double] (9,3)--(11,3);
    \draw [black,thick] (11,3)--(12,2);
    \draw [black,thick] (12,2)--(13,1);
    \plotpermnobox[black] {}{1,2,3,0,3,0,3,4,3,0,3,2,1}
\node at (1.325,1.675) {$1$};
    \node at (2.325,2.675) {$1$};
    \node at (4,3.25) {$1$};
    \node at (6,3.25) {$2$};
    \node at (7.325,3.675) {$1$};
    \node at (8.675,3.675) {$1$};
    \node at (10,3.25) {$3$};
    \node at (11.675,2.675) {$2$};
    \node at (12.675,1.675) {$1$};
    
\end{tikzpicture}
\end{center}
\vspace*{-2mm}
\caption{
   A super-augmented perfect matching of $[2n]$ with $n=6$,
   together with corresponding labeled Schr\"oder path
   (the label $\xi_i$ is written above the step $s_i$).
   The long level steps of type 2 are shown by double lines.
}
\label{fig:perfectmatchingwithpath}
\end{figure}

{\bf Step 1: Definition of the Schr\"oder path.}
Given a super-augmented perfect matching $\tau \in \scrm_{2n}^\star$,
we define a path $\omega = (\omega_0,\ldots,\omega_{2n})$
starting at $\omega_0 = (0,0)$ and ending at $\omega_{2n} = (2n,0)$,
with steps $s_1,\ldots,s_{2n}$ ending at locations $\omega_i = (i,h_i)$,
as follows:
\begin{itemize}
   \item If $i$ is a pure opener,
       then $s_i$ is a rise, so that $h_i = h_{i-1} + 1$.
   \item If $i$ is a pure closer,
       then $s_i$ is a fall, so that $h_i = h_{i-1} - 1$.
   \item If $(i,i+1)$ is a wiggly pair,
       then $s_i$ is a long level step of type~1
       (and $s_{i+1}$ is undefined).
       In this case the height $h_i$ is undefined,
       but we have $h_{i+1} = h_{i-1}$.
   \item If $(i,i+1)$ is a dashed pair,
       then $s_i$ is a long level step of type~2
       (and $s_{i+1}$ is undefined).
       In this case the height $h_i$ is undefined,
       but we have $h_{i+1} = h_{i-1}$.
\end{itemize}
(See Figure~\ref{fig:perfectmatchingwithpath} for an example.)
The interpretation of the heights $h_i$ is almost immediate
from this definition:

\begin{lemma}
   \label{lemma.heights.KZ}
For $i \in \{0,\ldots,2n\}$,
whenever the height $h_i$ is defined
it equals the number of arches that are ``started but unfinished''
after stage $i$
(or equivalently, before stage $i+1$), i.e.
\be
   h_i  \;=\;  \#\{ j \le i < k \colon\: j \sim_{\pi(\tau)} k \}
   \;.
\ee
\end{lemma}

\noindent
In particular, it follows that $\omega$ is indeed a Schr\"oder path,
i.e.\ all the heights $h_i$ are nonnegative (when they are defined)
and $h_{2n} = 0$.

\bigskip

\pagebreak
{\bf Step 2: Definition of the labels $\bm{\xi_i}$.}
\begin{itemize}
   \item If $i$ is a pure opener,
we set $\xi_i = 1$ as required by \reff{def.abc.KZ}.
   \item If $i$ is a pure closer,
we look at the $h_{i-1}$ arches
that are ``started but unfinished'' after stage $i-1$
(note that we must have $h_{i-1} \ge 1$);
let the openers of these arches be $x_1 < x_2 < \ldots < x_{h_{i-1}}$.
Then the vertex $i$ is paired with precisely one of these openers;
if it is $x_j$, we set $\xi_i = j$.
Obviously $1 \le \xi_i \le h_{i-1}$ as required by \reff{def.abc.KZ}.
   \item If $(i,i+1)$ is a wiggly pair,
then we look at the $h_{i-1}$ arches
that are ``started but unfinished'' after stage $i-1$,
and we define $\xi_i$ exactly as in the preceding case.
Again $1 \le \xi_i \le h_{i-1}$.
   \item If $(i,i+1)$ is a dashed pair,
then there are two possibilities [see Figure~\ref{fig:augmented}(b)]:
the opener $i$ and the closer $i+1$ could belong to different arches
(which necessarily cross), or they could belong to the same arch.
In the former case we look at the $h_{i-1}$ arches
that are ``started but unfinished'' after stage $i-1$,
and we define $\xi_i$ as before
(except that it is now vertex $i+1$ rather than $i$ that is paired with
 one of these openers).
In the latter case we set $\xi_i = h_{i-1} + 1$.
Obviously $1 \le \xi_i \le h_{i-1} + 1$ as required by \reff{def.abc.KZ}.
\end{itemize}
See again Figure~\ref{fig:perfectmatchingwithpath}.

\bigskip

{\bf Step 3: Proof of bijection.}
It is easy to describe the inverse map from
labeled Schr\"oder paths $(\omega,\xi)$ to super-augmented perfect matchings.
Successively for $i=1,\ldots,n$,
we use the 2-colored Schr\"oder path $\omega$
to read off the type associated to step $s_i$
(opener with no wiggly or dashed line, etc.).\footnote{
   It is a slight abuse of language here to call
   the 2-colored Schr\"oder path $\omega$,
   since by $\omega$ we mean here the sequence of {\em steps}\/ $s_i$,
   i.e.\ including also the type (1 or 2) of the long level steps.
   We trust that there will be no confusion.
}
And then, if step $s_i$ corresponds to anything other than
an opener with no wiggly or dashed line attached,
we use the label $\xi_i$ to decide to which opener
the vertex $i$ (or $i+1$) should be attached.

\bigskip

{\bf Step 4: Translation of the statistics.}

\begin{lemma}
   \label{lemma.statistics.KZ}
\quad\hfill
\vspace*{-1mm}
\begin{itemize}
   \item[(a)]  If $i$ is a pure opener, then
\be
   \qne(i,\pi)  \;=\;  h_{i-1}
   \;.
 \label{eq.lemma.statistics.KZ.opener}
\ee
   \item[(b)]  If $i$ is a pure closer, then
\begin{subeqnarray}
   \crr(i,\pi)   & = &   h_{i-1} - \xi_i  \\[2mm]
   \nee(i,\pi)   & = &   \xi_i - 1
 \label{eq.lemma.statistics.KZ.closer}
\end{subeqnarray}
   \item[(c)]  If $(i,i+1)$ is a wiggly pair, then also
\reff{eq.lemma.statistics.KZ.closer} holds.
   \item[(d)]  If $(i,i+1)$ is a dashed pair, then
\begin{subeqnarray}
   \crr(i+1,\pi)   & = &   h_{i-1} + 1 - \xi_i  \\[2mm]
   \nee(i+1,\pi)   & = &   \xi_i - 1
 \label{eq.lemma.statistics.KZ.openercloser}
\end{subeqnarray}
\end{itemize}
\end{lemma}

\proof
(a)  Each of the $h_{i-1}$ arches that are
``started but unfinished'' after stage $i-1$
will either cross or nest the arch that starts at $i$;
so this is an immediate consequence of the definition \reff{def.qne}.

(b,c) Look at the $h_{i-1}$ arches
that are ``started but unfinished'' after stage $i-1$,
and let the openers of these arches be $x_1 < x_2 < \ldots < x_{h_{i-1}}$;
by definition the vertex $i$ is paired with $x_{\xi_i}$.
Then each arch starting at a point $x_j$ with $j < \xi_i$
must nest with (and lie above) the arch from $x_{\xi_i}$ to $i$,
while each arch starting at a point $x_j$ with $j > \xi_i$
must cross the arch from $x_{\xi_i}$ to $i$.

(d) Let $x_1 < x_2 < \ldots < x_{h_{i-1}}$ be as before;
and let $x_{h_{i-1}+1} = i$.
By definition the vertex $i+1$ is paired with $x_{\xi_i}$.
Then the counting of nestings and crossings is exactly as in (b,c),
but with $h_{i-1}$ replaced by $h_{i-1} + 1$.
\qed

\bigskip

{\bf Step 5: Computation of the weights
   \reff{eq.thm.matchings.Ttype.final1.coeffs}.}
Using the bijection, we transfer the weights \reff{def.Mn}
from the super-augmented perfect matching $\tau$
to the labeled Schr\"oder path $(\omega,\xi)$
and then sum over $\xi$ to obtain the weight $W(\omega)$.
This weight is factorized over the individual steps $s_i$, as follows:
\begin{itemize}
   \item If $s_i$ is a rise starting at height $h_{i-1} = k$
      (so that $i$ is a pure opener),
      then from \reff{eq.lemma.statistics.KZ.opener} the weight is
\be
   a_k  \;=\;  \sfa_k   \;.
\ee
   \item If $s_i$ is a fall starting at height $h_{i-1} = k$
      (so that $i$ is a pure closer),
      then from \reff{eq.lemma.statistics.KZ.closer} the weight is
\be
   b_k
   \;=\;
   \sum_{\xi_i=1}^k \sfb_{k-\xi_i,\, \xi_i -1}
   \;=\;
   \sfb^\star_{k-1}
\ee
where $\sfb^\star_{k-1}$ was defined in \reff{def.bstar}.
   \item If $s_i$ is a long level step of type 1 at height $h_{i-1} = k$
      (so that $(i,i+1)$ is a wiggly pair and $k \ge 1$),
      then from \reff{eq.lemma.statistics.KZ.closer} the weight is
\be
   c^{(1)}_k
   \;=\;
   \sum_{\xi_i=1}^k \sff_{k-\xi_i,\, \xi_i -1}
   \;=\;
   \sff^\star_{k-1}
      \;.
\ee
   \item If $s_i$ is a long level step of type 2 at height $h_{i-1} = k$
      (so that $(i,i+1)$ is a dashed pair),
      then from \reff{eq.lemma.statistics.KZ.openercloser} the weight is
\be
   c^{(2)}_k
   \;=\;
   \sum_{\xi_i=1}^{k+1} \sfg_{k+1-\xi_i,\, \xi_i -1}
   \;=\;
   \sfg^\star_{k}
      \;.
\ee
\end{itemize}
Setting $\alpha_i = a_{i-1} b_i$ and
$\delta_i = c^{(1)}_{i-1} + c^{(2)}_{i-1}$
as instructed in \reff{eq.flajolet.schroder},
we obtain the weights \reff{eq.thm.matchings.Ttype.final1.coeffs}.
This completes the proof of Theorem~\ref{thm.matchings.Ttype.final1}.
\qed

\appendix
\section{Bijection between augmented perfect matchings and phylogenetic trees}
  \label{app.bijection}

We will now describe a bijection between augmented perfect matchings of $[2n]$
containing $\ell$ wiggly lines and phylogenetic trees with $n+1$
leaves and $n-\ell$ internal vertices. This bijection is illustrated in
Figure \ref{fig:bijection_perfect_matching_to_phylogenetic_tree}.

\usetikzlibrary{decorations.pathmorphing}
\tikzset{snake it/.style={decorate, decoration={snake,segment length=5,amplitude=1}}}


\newcommand{\plotpta}[4][] 
{ \fill[#1,radius={0.175*(#4)}] (#2,#3) circle; }


\newcommand{\plotpoints}[4][]  
{
  \foreach \y [count=\x] in {#3}
  {
    \ifnum0=\y {} \else {
      \plotpta[#1]{{(#2)*\x}}{\y}{#4}
    } \fi
  }
}

\newcommand{\WArc}[2] {\draw [black,very thick,snake it] (#2,1) arc (0:180:{0.5*(#2)-0.5*(#1)} and {0.33*(#2)-0.33*(#1)+0.2});}

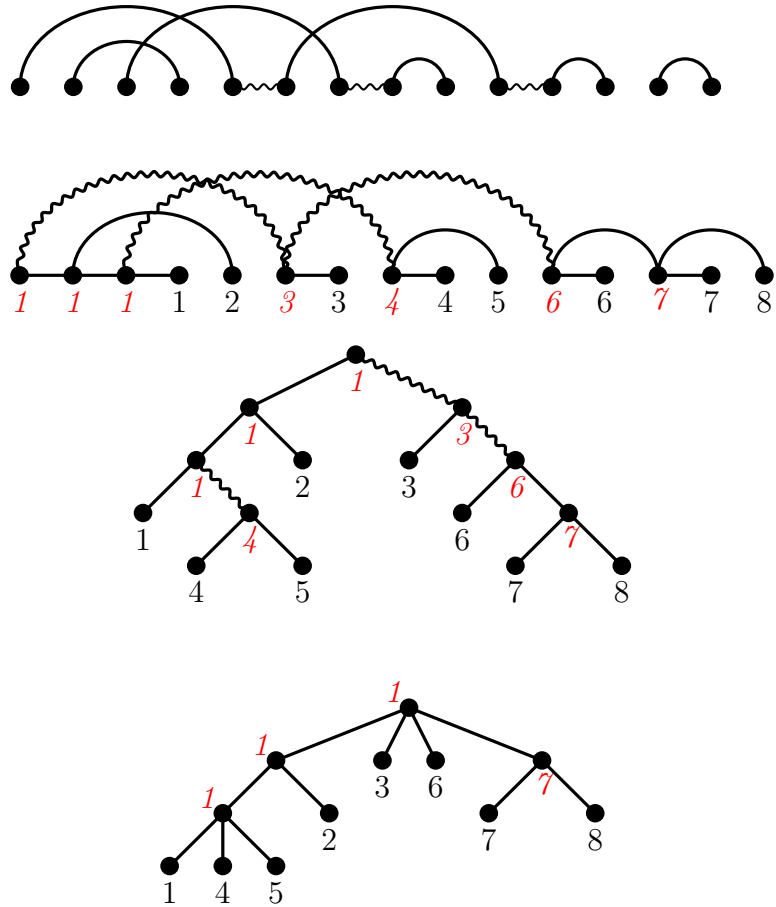
\begin{figure}[t]
$$
  \begin{tikzpicture}[scale=0.7,line join=round]
    \Arc{2}{4}
    \Arc{3}{7}
    \Arc{1}{5}
    \Arc{6}{10}
    \Arc{8}{9}
    \Arc{11}{12}
    \Arc{13}{14}
    \Wiggly{5}{6}
    \Wiggly{7}{8}
    \Wiggly{10}{11}
    \node at (14.5,0) {$~$};
    \plotpoints[black]{1}{1,1,1,1,1,1,1,1,1,1,1,1,1,1}{1}
  \end{tikzpicture}
$$
$$
    \hspace*{2.5mm}  
  \begin{tikzpicture}[scale=0.7,line join=round]
    \draw [black,very thick] (1,1)--(4,1);
    \draw [black,very thick] (6,1)--(7,1);
    \draw [black,very thick] (8,1)--(9,1);
    \draw [black,very thick] (11,1)--(12,1);
    \draw [black,very thick] (13,1)--(14,1);
    \node[red] at (1,0.5) {\it 1};
    \node[red] at (2,0.5) {\it 1};
    \node[red] at (3,0.5) {\it 1};
    \node at (4,0.5) {$1$};
    \node at (5,0.5) {$2$};
    \node[red] at (6,0.5) {\it 3};
    \node at (7,0.5) {$3$};
    \node[red] at (8,0.5) {\it 4};
    \node at (9,0.5) {$4$};
    \node at (10,0.5) {$5$};
    \node[red] at (11,0.5) {\it 6};
    \node at (12,0.5) {$6$};
    \node[red] at (13,0.5) {\it 7};
    \node at (14,0.5) {$7$};
    \node at (15,0.5) {$8$};
    \Arc{2}{5}
    \WArc{3}{8}
    \WArc{1}{6}
    \WArc{6}{11}
    \Arc{8}{10}
    \Arc{11}{13}
    \Arc{13}{15}
    \plotpoints[black]{1}{1,1,1,1,1,1,1,1,1,1,1,1,1,1,1}{1}
  \end{tikzpicture}
$$
    $$
      \begin{tikzpicture}[scale=0.7,line join=round]
    \node at (1, 1.5){$1$};
    \node at (4, 2.5){$2$};
    \node at (6, 2.5){$3$};
    \node at (2, 0.5){$4$};
    \node at (4, 0.5){$5$};
    \node at (7, 1.5){$6$};
    \node at (8, 0.5){$7$};
    \node at (10, 0.5){$8$};
    \node [red] at (2, 2.5){\it 1};
    \node [red] at (3, 3.5){\it 1};
    \node [red] at (5, 4.5){\it 1};
    \node [red] at (7, 3.5){\it 3};
    \node [red] at (8, 2.5){\it 6};
    \node [red] at (9, 1.5){\it 7};
    \node [red] at (3, 1.5){\it 4};
    \plotpt[black]{5}{5}
    \plotpt[black]{2}{3}
    \plotpt[black]{3}{4}
    \plotpt[black]{4}{3}
    \plotpt[black]{1}{2}
    \plotpt[black]{6}{3}
    \plotpt[black]{7}{4}
    \plotpt[black]{2}{1}
    \plotpt[black]{3}{2}
    \plotpt[black]{8}{3}
    \plotpt[black]{4}{1}
    \plotpt[black]{9}{2}
    \plotpt[black]{7}{2}
    \plotpt[black]{10}{1}
    \plotpt[black]{8}{1}
    \draw[black,very thick] (5, 5)--(3, 4);
    \draw[black,very thick] (3, 4)--(2, 3);
    \draw[black,very thick] (2, 3)--(1, 2);
    \draw[black,very thick] (3, 4)--(4, 3);
    \draw[black,very thick,snake it] (5, 5)--(7, 4);
    \draw[black,very thick] (7, 4)--(6, 3);
    \draw[black,very thick,snake it] (2, 3)--(3, 2);
    \draw[black,very thick] (3, 2)--(2, 1);
    \draw[black,very thick] (3, 2)--(4, 1);
    \draw[black,very thick,snake it] (7, 4)--(8, 3);
    \draw[black,very thick] (8, 3)--(7, 2);
    \draw[black,very thick] (8, 3)--(9, 2);
    \draw[black,very thick] (9, 2)--(8, 1);
    \draw[black,very thick] (9, 2)--(10, 1);

      \end{tikzpicture}
    $$

    $$
      \begin{tikzpicture}[scale=0.7,line join=round]
    \node at (1, 1.5){$1$};
    \node at (4, 2.5){$2$};
    \node at (5, 3.5){$3$};
    \node at (2, 1.5){$4$};
    \node at (3, 1.5){$5$};
    \node at (6, 3.5){$6$};
    \node at (7, 2.5){$7$};
    \node at (9, 2.5){$8$};
    \node[red] at (1.7, 3.3){\it 1};
    \node[red] at (2.7, 4.3){\it 1};
    \node[red] at (5.2, 5.3){\it 1};
    \node[red] at (8, 3.5){\it 7};
    \plotpt[black]{5.5}{5}
    \plotpt[black]{2}{3}
    \plotpt[black]{3}{4}
    \plotpt[black]{4}{3}
    \plotpt[black]{1}{2}
    \plotpt[black]{5}{4}
    \plotpt[black]{2}{2}
    \plotpt[black]{3}{2}
    \plotpt[black]{8}{4}
    \plotpt[black]{6}{4}
    \plotpt[black]{9}{3}
    \plotpt[black]{7}{3}
    \draw[black,very thick] (5.5, 5)--(3, 4);
    \draw[black,very thick] (3, 4)--(2, 3);
    \draw[black,very thick] (2, 3)--(1, 2);
    \draw[black,very thick] (3, 4)--(4, 3);
    \draw[black,very thick] (5.5, 5)--(5, 4);
    \draw[black,very thick] (2, 3)--(2, 2);
    \draw[black,very thick] (2, 3)--(3, 2);
    \draw[black,very thick] (5.5, 5)--(6, 4);
    \draw[black,very thick] (5.5, 5)--(8, 4);
    \draw[black,very thick] (8, 4)--(7, 3);
    \draw[black,very thick] (8, 4)--(9, 3);

      \end{tikzpicture}
    $$
\caption{
   An augmented perfect matching $P$ (top row) and the corresponding
   phylogenetic tree $T$ (bottom row), along with the intermediate steps
   in the bijection: the arch system $P'$ (second row)
   and the planar rooted binary tree $T_2$ (third row).
   In the forward direction of the bijection, only the black parts are used.
   The extra labels used in the reverse bijection are shown in red italic.
}
  \label{fig:bijection_perfect_matching_to_phylogenetic_tree}
\end{figure}

Given an augmented perfect matching $P$ of $[2n]$,
we construct an arch system $P'$ on $[2n+1]$ as follows:
For each arch $(i,j)$ in $P$, we draw the arch $(i,j+1)$ in $P'$
as well as a horizontal line $(i,i+1)$.
If $j$ and $j+1$ are joined by a wiggly line in $P$,
we make the arch $(i,j+1)$ wiggly in $P'$.
Now each vertex $i\geq2$ in $P'$
is joined to exactly one vertex on its left:
if $i-1$ is an opener in $P$,
then $i$ is joined to $i-1$ by a horizontal line;
whereas if $i-1$ is a closer in $P$,
then $i$ is the closer of an arch in $P'$.
Hence the graph formed by $P'$ is a tree.
Moreover, each opener in $P'$ is joined to exactly two vertices on its right:
one by an arch and one by a horizontal edge.
And each non-opener in $P'$ is not joined to any vertices on its right.

We now label the vertices in $P'$ that are not openers
by using the numbers $1,\ldots,n+1$ in order from left to right
(shown in black roman font in
 Figure~\ref{fig:bijection_perfect_matching_to_phylogenetic_tree}).
Each labeled vertex in $P'$ is not joined to any vertices on its right.
Hence we can define a planar rooted binary tree $T_{2}$
that is isomorphic as a labeled graph to $P'$,
such that horizontal edges in $P'$ correspond to left edges in $T_{2}$,
while arches in $P'$ correspond to right edges of $T_{2}$
(which are wiggly whenever the corresponding arch is wiggly).
Since wiggly edges in the perfect matching join closers to openers,
the child of a wiggly edge in $T_{2}$ cannot be a leaf.

Finally, to construct the phylogenetic tree $T$,
we simply contract all wiggly edges of $T_{2}$.
Since the child of a wiggly edge in $T_{2}$ cannot be a leaf,
each internal vertex in $T$ has at least two children,
so $T$ is indeed a phylogenetic tree.

To show that this is a bijection, we will describe the reverse transformation.
Given a phylogenetic tree $T$, we start by labeling each internal vertex of $T$
with the minimum label amongst its descendants
(shown in red italic in
 Figure~\ref{fig:bijection_perfect_matching_to_phylogenetic_tree}).
We then order the children of each vertex from left to right in
increasing order of their label. It is easy to see that the tree
$T$ constructed from an augmented perfect matching will always be
drawn in this way. Now we will describe how to construct the tree
$T_{2}$ from $T$. For each vertex $v$ with degree at least $3$, let
$c_{1},c_{2},\ldots,c_{k}$ be the children of $v$ in order from left
to right (that is, in increasing order of labels). Then we split $v$
into a sequence of vertices $v_{1},\ldots,v_{k-1}$ of out-degree~$2$
so that the left child of $v_{i}$ is $c_{i}$ and the right child of
$v_{i}$ is $v_{i+1}$ if $i<k-1$, while the right child of $v_{k-1}$
is $c_{k}$. Moreover, each edge joining $v_{i}$ to its right child
$v_{i+1}$ is wiggly. We then label each new vertex with the same label
as its left child. To construct $P'$ from $T_{2}$, we just have to
order the vertices from left to right. If $u$ and $v$ are vertices of
$T_{2}$ such that $u$ has a lower label than $v$, then we say $u<v$,
and if $v$ is the left child of $u$ (so they have the same label)
then we also say $u<v$. We construct $P'$ by placing the vertices in
increasing order according to $<$. Finally $P$ is constructed from
$P'$ by creating an arch $(i,j)$ for each arch $(i,j+1)$ in $P'$,
and a wiggly line $(j,j+1)$ for each wiggly arch $(i,j+1)$ in $P'$.

Since the transformations between $P$ and $T$ are inverses, each transformation is a bijection.

\section{Recurrence for polynomials defined by the general linear T-fraction}
   \label{app.recurrence}

Consider the T-fraction with coefficients $\alpha_i = x + (i-1)u$
and $\delta_i = z + (i-1)w$,
and let $W_n(x,u,z,w)$ be the polynomials that it generates:
\begin{subeqnarray}
   f(t;x,u,z,w)
   & \eqdef &
   \sum_{n=0}^\infty W_n(x,u,z,w) \, t^n
             \\[2mm]
   & \eqdef &
   \cfrac{1}{1 - zt - \cfrac{xt}{1 - (z+w)t - \cfrac{(x+u)t}{1 - (z+2w)t - \cfrac{(x+2u)t}{1- \cdots}}}}
   \;.
   \qquad
  \label{eq.linear.T-fraction}
\end{subeqnarray}
We will prove:

\begin{proposition}
   \label{prop.recurrence}
The ordinary generating function $f(t;x,u,z,w)$ satisfies the
nonlinear partial differential equation
\be
   f  \;=\;  1 \,+\, (u+z)t f \,+\, ut^2 f_t \,+\,
                     (u+w)t (u f_u + x f_x) \,+\, (x-u)t f^2
   \;.
\ee
Equivalently, the polynomials $W_n(x,u,z,w)$ satisfy the
nonlinear differential recurrence
\be
   W_n  \;=\;  \delta_{n0}  \,+\,
               (z+nu) W_{n-1} \,+\,
               (u+w) \Bigl( u {\partial W_{n-1} \over \partial u}
                            +
                            x {\partial W_{n-1} \over \partial x}
                     \Bigr)
               \,+\,
               (x-u) \sum_{j=0}^{n-1} W_j W_{n-1-j}
\ee
where $W_{-1} \eqdef 0$.
\end{proposition}

In particular, when we restrict to $u=x$,
we obtain a {\em linear}\/ partial differential equation
and a linear differential recurrence:

\begin{corollary}
   \label{cor.recurrence.1}
The ordinary generating function $g(t;x,z,w) \eqdef f(t;x,x,z,w)$
satisfies the linear partial differential equation
\be
   g  \;=\;  1 \,+\, (x+z)t g \,+\, ut^2 g_t \,+\,
                     (u+w)t (u g_u + x g_x)
   \;.
\ee
Equivalently, the polynomials $P_n(x,z,w) \eqdef W_n(x,x,z,w)$
satisfy the linear differential recurrence
\be
   P_n  \;=\;  \delta_{n0}  \,+\, (z+nx) P_{n-1} \,+\,
               x(x+w) {\partial P_{n-1} \over \partial x}
 \label{eq.cor.recurrence.1.Pn}
\ee
where $P_{-1} \eqdef 0$.
Equivalently, if we write
$W_n(x,x,z,w) = \sum\limits_{k=0}^n W_{n,k}(z,w) \, x^k$,
then the polynomials $W_{n,k}(z,w)$ satisfy the linear recurrence
\be
   W_{n,k}  \;=\;  (n+k-1) \, W_{n-1,k-1}  \:+\: (z+kw) \, W_{n-1,k}
   \quad\hbox{for $n \ge 1$}
 \label{eq.ward.recurrence.zw}
\ee
with initial condition $W_{0,k} = \delta_{k0}$.
\end{corollary}

\noindent
When $z=0$ and $w=1$, the recurrence \reff{eq.ward.recurrence.zw}
reduces to \reff{eq.ward.recurrence}.

On the other hand, when we restrict to $w = -u$,
the ordinary generating function satisfies
a nonlinear {\em ordinary}\/ differential equation of Riccati type:

\begin{corollary}
   \label{cor.recurrence.2}
The ordinary generating function $h(t;x,u,z) \eqdef f(t;x,u,z,-u)$
satisfies the Riccati equation
\be
   h  \;=\;  1 \,+\, (u+z)t h \,+\, ut^2 h' \,+\, (x-u)t h^2
\ee
where ${}'$ denotes $\partial/\partial t$.
Equivalently, the polynomials $Q_n(x,u,z) \eqdef W_n(x,u,z,-u)$
satisfy the nonlinear recurrence
\be
   Q_n  \;=\;  \delta_{n0}  \,+\,
               (z+nu) Q_{n-1} \,+\,
               (x-u) \sum_{j=0}^{n-1} Q_j Q_{n-1-j}
\ee
where $Q_{-1} \eqdef 0$.
\end{corollary}
   
And finally, when we restrict to $u=0$,
the ordinary generating function satisfies
a nonlinear ordinary differential equation of Riccati type,
but in the variable $x$ rather than $t$:

\begin{corollary}
   \label{cor.recurrence.3}
The ordinary generating function $H(x;t,z,w) \eqdef f(t;x,0,z,w)$
satisfies the Riccati equation
\be
   H  \;=\;  1 \,+\, zt H \,+\, wxt H_x \,+\, xt H^2
   \;.
\ee
\end{corollary}

\proofof{Proposition~\ref{prop.recurrence}}
A combinatorial interpretation of the polynomials $W_n(x,u,z,w)$
was given in Theorem~\ref{thm1.2}:
namely, $W_n(x,u,z,w)$ is the generating polynomial for
super-augmented perfect matchings of $[2n]$ in which
each pure closer with crossing number 0 gets a weight~$x$,
each pure closer with crossing number $\ge 1$ gets a weight~$u$,
each dashed line for which the two endpoints belong to the same arch
gets a weight~$z$,
each other dashed line gets a weight~$w''$,
and each wiggly line gets a weight~$w'$,
whenever $w = w' + w''$.
Here we use the interpretation with $w' = 0$,
so that wiggly lines are forbidden.

The contribution to $f$ from the case with no arches ($n=0$) is clearly~1.
Otherwise, let the arch whose opener is vertex~1 be called $\alpha$,
and call the closer of this arch $c_{\alpha}$.
Let $A'$ be the perfect matching that remains when $\alpha$,
its incident vertices and their incident dashed lines are removed,
and let $t^{a+b+c+d}x^au^bz^cw^d$ be the contribution of $A'$ to $f$. We will proceed by considering five cases for the type of $\alpha$,
as illustrated in Figure~\ref{fig:5cases}.

\newcommand{\fillArch}[2] {\fill [gray!40!white] (#2,1) arc (0:180:{0.5*(#2)-0.5*(#1)}); \draw [gray!40!white,line width=6] (#2,1)--(#1,1);}

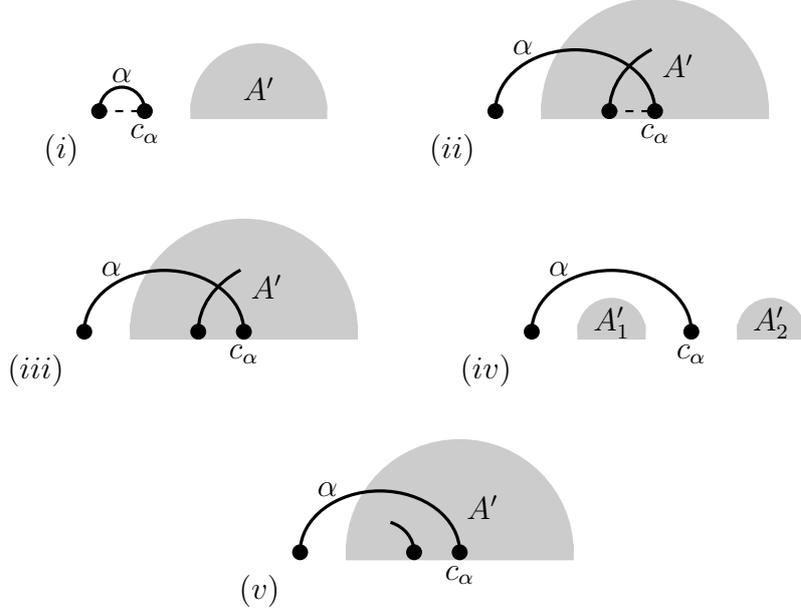
\begin{figure}[t]
$$
  (i)\begin{tikzpicture}[scale=0.6,line join=round]
    \fillArch{3}{6}
  	\Arc{1}{2}
    \node at (2,0.5){$c_{\alpha}$};
    \node at (1.5,1.8){$\alpha$};
  	\Dashed{1}{2}
    \node at (4.5, 1.5){$A'$};
    \node at (0.9,0.2){$~$};
    \plotpoints[black]{1}{1,1}{1}
    
  \end{tikzpicture}
\qquad\quad
  (ii)\begin{tikzpicture}[scale=0.6,line join=round]
    \fillArch{2}{7}
  	\Arc{1}{4.5}
    \node at (4.5,0.5){$c_{\alpha}$};
    \node at (1.6,2.4){$\alpha$};
  	\Dashed{3.5}{4.5}
    \node at (5, 2){$A'$};
    \plotpoints[black]{0.5}{0,1,0,0,0,0,1,0,1}{1}
    \PartArc{-1}{3.5}{0.5}{0.35}
    \node at (0.9,0.2){$~$};
    
  \end{tikzpicture}
$$
$$
  (iii)\begin{tikzpicture}[scale=0.6,line join=round]
    \fillArch{2}{7}
  	\Arc{1}{4.5}
    \node at (4.5,0.5){$c_{\alpha}$};
    \node at (1.6,2.4){$\alpha$};
    \node at (5, 2){$A'$};
    \plotpoints[black]{0.5}{0,1,0,0,0,0,1,0,1}{1}
    \PartArc{-1}{3.5}{0.5}{0.35}
    \node at (0.9,0.2){$~$};
    \node at (4,4){$~$};
    
  \end{tikzpicture}
\qquad\quad
  (iv)\begin{tikzpicture}[scale=0.6,line join=round]
    \fillArch{2}{3.5}
    \fillArch{5.5}{7}
  	\Arc{1}{4.5}
    \node at (4.5,0.5){$c_{\alpha}$};
    \node at (1.6,2.4){$\alpha$};
    \node at (2.75, 1.2){$A_{1}'$};
    \node at (6.25, 1.2){$A_{2}'$};
    \plotpoints[black]{0.5}{0,1,0,0,0,0,0,0,1}{1}
    \node at (0.9,0.2){$~$};
    
  \end{tikzpicture}
$$
$$
  (v)\begin{tikzpicture}[scale=0.6,line join=round]
    \fillArch{2}{7}
  	\Arc{1}{4.5}
    \node at (4.5,0.5){$c_{\alpha}$};
    \node at (1.6,2.4){$\alpha$};
    \node at (5, 2){$A'$};
    \plotpoints[black]{0.5}{0,1,0,0,0,0,1,0,1}{1}
    \PartArc{2}{3.5}{0}{0.2}
    \node at (0.9,0.2){$~$};
    \node at (4,4){$~$};
    
  \end{tikzpicture}
$$
\caption{
   The five cases in the proof of Proposition~\ref{prop.recurrence}.
}
\label{fig:5cases}
\end{figure}

Case~$(i)$ is when vertex~1 is incident on a dashed line.
Then vertex~2 is a closer, so it must close $\alpha$, i.e.\ $c_\alpha = 2$.
Then $A'$ can be any dash-augmented arch system,
so the contribution from this case is $ztf$.

In the remaining four cases, vertex~1 is a pure opener.
Case~$(ii)$ is when $c_{\alpha}$ is incident on a dashed line.
Then $\alpha$ contributes the weight $wt$.
The closer $c_{\alpha}$ can be placed immediately after
any pure opener in $A'$, of which there are $a+b$.
Hence the contribution to $f$ from arch systems corresponding to
the smaller arch system $A'$ is $(a+b)t^{a+b+c+d+1}x^au^bz^cw^{d+1}$.
Summing this over all possible arch systems $A'$
yields the contribution $wt(uf_{u}+xf_{x})$ from this case.

In the remaining three cases, $c_{\alpha}$ is a pure closer.
Case~$(iii)$ is when $c_{\alpha}$ immediately follows the opener
of a different arch
(then $c_{\alpha}$ necessarily has crossing number $\ge 1$).
This is identical to Case~$(ii)$ except that
we do not attach a dashed line to $c_{\alpha}$,
so $\alpha$ contributes the weight $ut$ instead of $wt$.
Hence the contribution from  this case is $ut(uf_{u}+xf_{x})$.

Case~$(iv)$ is when $\alpha$ has crossing number~0.
Then the arch system $A'$ separates into two sections $A_{1}'$ and $A_{2}'$,
as shown in Figure \ref{fig:5cases},
where $A_{1}'$ is the section contained under $\alpha$
and $A_{2}'$ is the section which follows $c_{\alpha}$.
The arch $\alpha$ contributes the weight $xt$,
so the contribution from this case is $xtf^2$.

Case~$(v)$ is when $\alpha$ has crossing number $\ge 1$
and $c_{\alpha}$ immediately follows another closer.
If we remove the condition that $\alpha$ has crossing number $\ge 1$,
then $c_{\alpha}$ can follow any closer in $A'$,
so there are $a+b+c+d$ possible positions for $c_{\alpha}$.
Ignoring the weight contributed by $\alpha$, these are counted by $tf_{t}$.
Now we subtract $f^2$ from this to remove the cases where
$\alpha$ has crossing number $0$ (namely, Case~$(iv)$);
however, we add back $f$ to account for the case in which $c_{\alpha}$
is vertex~2.
Multiplying all this by the weight $ut$ of $\alpha$
yields the contribution $ut(tf_{t}-f^2+f)$ for this case.

Adding the contributions from all five cases yields the desired result.
\qed

\bigskip

Finally, let us use Corollary~\ref{cor.recurrence.1} to prove the formula
\reff{eq.multiward.inverse.Wnxuzw.u=x}/\reff{eq.multiward.inverse.Wnxuzw.u=x.F}
corresponding to the case $u=x$:

\proofof{\reff{eq.multiward.inverse.Wnxuzw.u=x}/\reff{eq.multiward.inverse.Wnxuzw.u=x.F}}
Multiply \reff{eq.cor.recurrence.1.Pn} by $t^n/n!$ and sum over $n \ge 0$:
this shows that the exponential generating function
\be
   \scrw(t;x,z,w)
   \;\eqdef\;
   \sum_{n=0}^\infty W_n(x,x,z,w) \, {t^{n+1} \over (n+1)!}
 \label{def.scrw.bis}
\ee
[cf.\ \reff{def.scrw}] satisfies the linear partial differential equation
\be
   \scrw_t  \;=\;  1 \,+\, z\scrw \,+\, xt \scrw_t \,+\, x(x+w) \scrw_x
   \;.
 \label{eq.pde.scrw}
\ee
Now let $F(t;x,z,w)$ be the compositional inverse of $\scrw(t;x,z,w)$
with respect to $t$, which satisfies
\be
   \scrw\bigl( F(t;x,z,w);\, x,z,w \bigr)  \;=\; t
   \;.
 \label{eq.scrw.F}
\ee
Differentiation of \reff{eq.scrw.F} yields
$\scrw_t = 1/F_t$ and $\scrw_x = -F_x/F_t$
(where of course $\scrw$ is evaluated at $t \leftarrow F(t;x,z,w)$).
Evaluating \reff{eq.pde.scrw} at $t \leftarrow F(t;x,z,w)$
then shows that $F$ satisfies the linear partial differential equation
\be
   F_t  \;=\;  1 \,-\, xF \,-\, zt F_t \,+\, x(x+w) F_x
   \;.
 \label{eq.pde.F}
\ee
The function $F(t;x,z,w)$ is uniquely determined
by this partial differential equation
together with the initial condition $F(0;x,z,w) = 0$.
And it is straightforward to verify that the expression
\reff{eq.multiward.inverse.Wnxuzw.u=x.F}
indeed satisfies \reff{eq.pde.F}.
\qed

%

\clearpage
\section*{Acknowledgments}

We wish to thank Bishal Deb, Mathias P\'etr\'eolle and Jiang Zeng
for helpful conversations.

This work has benefited greatly from the existence of
the On-Line Encyclopedia of Integer Sequences \cite{OEIS}.
We warmly thank Neil Sloane for founding this indispensable resource,
and the hundreds of volunteers for helping to maintain and expand it.

This research was supported in part by
Engineering and Physical Sciences Research Council grant EP/N025636/1.
The first author also received funding from the European Research
Council (ERC) under the European Union's Horizon 2020 research and
innovation programme under the Grant Agreement No.~759702.


\begin{thebibliography}{99}

\bibitem{Barbero_14}  J.F. Barbero G., J. Salas and E.J.S. Villase\~nor,
   Bivariate generating functions for a class of linear recurrences:
   General structure, J. Combin. Theory A {\bf 125}, 146--165 (2014).

\bibitem{Barbero_15}  J.F. Barbero G., J. Salas and E.J.S. Villase\~nor,
   Generalized Stirling permutations and forests: higher-order Eulerian
   and Ward numbers, Electron. J. Combin. {\bf 22}, no.~3, paper 3.37 (2015).

\bibitem{Barry_09}  P. Barry, Continued fractions and transformations of
   integer sequences, J. Integer Seq. {\bf 12}, article 09.7.6 (2009).


\bibitem{Blitvic_12}  N. Blitvi\'c, The $(q,t)$-Gaussian process,
   J. Funct. Anal. {\bf 263}, 3270--3305 (2012).

\bibitem{Callan_09}  D. Callan, A combinatorial survey of identities
   for the double factorial,
   preprint (2009), arXiv:0906.1317 [math.CO] at arXiv.org.

\bibitem{Callan_15}  D. Callan, T. Mansour and M. Shattuck,
   Some identities for derangement and Ward number sequences and
   related bijections,
   Pure Math. Appl. (PU.M.A.) {\bf 25}, 132--143 (2015).

\bibitem{Carlitz_71}  L. Carlitz, Note on the numbers of Jordan and Ward,
   Duke Math. J. {\bf 38}, 783--790 (1971).


\bibitem{Clark_99}  L. Clark, Asymptotic normality of the Ward numbers,
   Discrete Math. {\bf 203}, 41--48 (1999).

\bibitem{Comtet_70}  L. Comtet, Sur le quatri\`eme probl\`eme et les nombres
   de Schr\"oder,
   Comptes Rendus Acad. Sci. Paris S\'er. A-B {\bf 271}, A913--A916 (1970).

\bibitem{Comtet_74}  L. Comtet, {\em Advanced Combinatorics:
   The Art of Finite and Infinite Expansions}\/
   (Reidel, Dordrecht--Boston, 1974).
   [French original:  {\em Analyse Combinatoire}\/, tomes~I et II,
   Presses Universitaires de France, Paris, 1970.]

\bibitem{Diaconis_98}  P.W. Diaconis and S.P. Holmes,
   Matchings and phylogenetic trees,
   Proc. Natl. Acad. Sci. USA {\bf 95}, 14600--14602 (1998).

\bibitem{Drake_07}  B. Drake, I.M. Gessel and G. Xin,
   Three proofs and a generalization of the Goulden--Litsyn--Shevelev
   conjecture on a sequence arising in algebraic geometry,
   J. Integer Seq. {\bf 10}, article 07.3.7 (2007).


\bibitem{Enestrom_13}  G. Enestr\"om, 
   {\em Die Schriften Eulers chronologisch nach den Jahren geordnet,
    in denen sie verfa{\ss}t worden sind}\/,
   Jahresbericht der Deutschen Mathematiker-Vereinigung
   (Teubner, Leipzig, 1913).

\bibitem{Erdos_89}  P.L. Erd\H{o}s and L.A. Sz\'ekely,
   Applications of antilexicographic order.\ I.~An enumerative theory of trees,
   Adv. Appl. Math. {\bf 10}, 488--496 (1989).
   
\bibitem{Euler_1748}  L. Euler, {\em Introductio in Analysin Infinitorum}\/,
   tomus primus (Bousquet, Lausanne, 1748).
   English translation:  L. Euler, {\em Introduction to Analysis of the
    Infinite, Book~I}\/, translated from the Latin and with an introduction
    by John D.~Blanton (Springer-Verlag, New York, 1988).
   [Latin original and French and German translations available at
    \url{http://eulerarchive.maa.org/pages/E101.html}]

\bibitem{Euler_1760}  L. Euler, De seriebus divergentibus,
   Novi Commentarii Academiae Scientiarum Petropolitanae {\bf 5}, 205--237
   (1760).
   [Latin original and English and German translations available at
    \url{http://eulerarchive.maa.org/pages/E247.html}]

\bibitem{Euler_1788}  L. Euler, De transformatione seriei divergentis
   $1 - mx + m(m+n)x^2 - m(m+n)(m+2n)x^3 + \hbox{etc.}$
   in fractionem continuam,
   Nova Acta Academiae Scientarum Imperialis Petropolitanae
     {\bf 2}, 36--45 (1788).
   [Latin original and English and German translations available at
    \url{http://eulerarchive.maa.org/pages/E616.html}]

\bibitem{Flajolet_80}  P. Flajolet, Combinatorial aspects of continued
   fractions,  Discrete Math. {\bf 32}, 125--161 (1980).


\bibitem{Flajolet_09}  P. Flajolet and R. Sedgewick,
  {\em Analytic Combinatorics}\/ (Cambridge University Press, Cambridge, 2009).

\bibitem{Flajolet_94}  P. Flajolet, P. Zimmerman and B. Van Cutsem,
   A calculus for the random generation of labelled combinatorial structures,
   Theoret. Comput. Sci. {\bf 132}, 1--35 (1994).


\bibitem{Fusy_15}  E. Fusy and E. Guitter,
   Comparing two statistical ensembles of quadrangulations:
   A continued fraction approach,
   Ann. Inst. Henri Poincar\'e D {\bf 4}, 125--176 (2017).

\bibitem{Gaiffi_15}  G. Gaiffi,
   Nested sets, set partitions and Kirkman--Cayley dissection numbers,
   European J. Combin. {\bf 43}, 279--288 (2015).

\bibitem{Gessel_78}  I. Gessel and R.P. Stanley, Stirling polynomials,
   J. Combin. Theory A {\bf 24}, 24--33 (1978).

\bibitem{Gladkovskii_13b}  S.N. Gladkovskii, 1~May 2013, contribution to
   \cite[sequence~A000311]{OEIS}.

\bibitem{Graham_94}  R.L. Graham, D.E. Knuth and O. Patashnik,
   {\em Concrete Mathematics: A Foundation for Computer Science}\/,
   2nd ed.~(Addison-Wesley, Reading, Mass., 1994).

\bibitem{Haiman_89}  M. Haiman and W. Schmitt,
   Incidence algebra antipodes and Lagrange inversion in one and several
   variables, J. Combin. Theory A {\bf 50}, 172--185 (1989).


\bibitem{Josuat-Verges_18}  M. Josuat-Verg\`es,
   A $q$-analog of Schl\"afli and Gould identities on Stirling numbers,
   Ramanujan J. {\bf 46}, 483--507 (2018).


\bibitem{Kasraoui_06}  A. Kasraoui and J. Zeng, Distribution of crossings,
   nestings and alignments of two edges in matchings and partitions,
   Electron. J. Combin. {\bf 13}, \#R33 (2006).

\bibitem{Leclerc_85}  B. Leclerc, Les hi\'erarchies de parties et leur
   demi-treillis,  Math. Sci. Humaines {\bf 89}, 5--34, 67 (1985);
   errata {\bf 92}, 40 (1985).

\bibitem{OEIS}  The On-Line Encyclopedia of Integer Sequences,
   published electronically at \url{http://oeis.org}

\bibitem{Oste_15}  R. Oste and J. Van der Jeugt,
   Motzkin paths, Motzkin polynomials and recurrence relations,
   Electron. J. Combin. {\bf 22}, no.~2, \#P2.8 (2015).

\bibitem{latpath_SRTR}  M. P\'etr\'eolle, A.D. Sokal and B.-X. Zhu,
   Lattice paths and branched continued fractions:
       An infinite sequence of generalizations
       of the Stieltjes--Rogers and Thron--Rogers polynomials,
       with coefficientwise Hankel-total positivity,
   preprint (2018), arXiv:1807.03271 [math.CO] at arXiv.org.

\bibitem{Prodinger_18}  H. Prodinger, A bijection between phylogenetic trees
   and plane oriented recursive trees,
   Rend. Istit. Mat. Univ. Trieste {\bf 50}, 133--137 (2018).


\bibitem{Riordan_68}  J. Riordan, {\em Combinatorial Identities}\/
    (Wiley, New York, 1968).
    [Reprinted with corrections by Robert E.~Krieger Publishing Co.,
     Huntington NY, 1979.]


\bibitem{Riordan_76}  J. Riordan,
   The blossoming of Schr\"oder's fourth problem,
   Acta Math. {\bf 137}, 1--16 (1976).

\bibitem{Roblet_96}  E. Roblet and X.G. Viennot,
   Th\'eorie combinatoire des T-fractions et approximants de Pad\'e
   en deux points,
   Discrete Math. {\bf 153}, 271--288 (1996).

\bibitem{Schroder_1870}  E. Schr\"oder, Vier kombinatorische Probleme,
   Z. f\"ur Math. Physik {\bf 15}, 361--376 (1870).

\bibitem{Sokal_flajolet}  A.D. Sokal, Coefficientwise total positivity
   (via continued fractions) for some Hankel matrices of combinatorial
   polynomials, talk at the S\'eminaire de Combinatoire Philippe Flajolet,
   Institut Henri Poincar\'e, Paris, 5 June 2014;
   transparencies available at
   \url{http://semflajolet.math.cnrs.fr/index.php/Main/2013-2014}

\bibitem{Sokal_totalpos}  A.D. Sokal, Coefficientwise total positivity
   (via continued fractions) for some Hankel matrices of combinatorial
   polynomials, in preparation.

\bibitem{Sokal-Zeng_masterpoly}  A.D. Sokal and J. Zeng,
   Some multivariate master polynomials for permutations, set partitions,
   and perfect matchings, and their continued fractions,
   in preparation.


\bibitem{Stanley_99}  R.P. Stanley, {\em Enumerative Combinatorics}\/,
      vol.~2 (Cambridge University Press, Cambridge--New York, 1999).

\bibitem{Steel_14}  M. Steel, Tracing evolutionary links between species,
   Amer. Math. Monthly {\bf 121}, 771--792 (2014).

\bibitem{Viennot_83}  G. Viennot, Une th\'eorie combinatoire des polyn\^omes
   orthogonaux g\'en\'eraux, Notes de conf\'erences donn\'ees
   \`a l'Universit\'e du Qu\'ebec \`a Montr\'eal,
   septembre-octobre 1983.
   Available on-line at
   \url{http://www.xavierviennot.org/xavier/polynomes_orthogonaux.html}

\bibitem{Wall_48}  H.S. Wall, {\em Analytic Theory of Continued Fractions}\/
   (Van Nostrand, New York, 1948).

\bibitem{Ward_34}  M. Ward, The representation of Stirling's numbers and
   Stirling's polynomials as sums of factorials,
   Amer. J. Math. {\bf 56}, 87--95 (1934).

\end{thebibliography}
\end{document}